\documentclass[12pt]{article}
\usepackage{color}
\usepackage{amscd,amsfonts,amssymb,amscd,amsmath,latexsym}
\usepackage[hypertex]{hyperref}
\usepackage[all]{xy}
\usepackage{mathrsfs}
\title{String structures and trivialisations of a Pfaffian line bundle}
\author{Ulrich Bunke\thanks{NWF I - Mathematik,
Universit{\"a}t Regensburg,
93040 Regensburg,
GERMANY, ulrich.bunke@mathematik.uni-regensburg.de} }

\newtheorem{theorem}{Theorem}[section] 
\newtheorem{prop}[theorem]{Proposition}
\newtheorem{lem}[theorem]{Lemma}
\newtheorem{ddd}[theorem]{Definition}
\newtheorem{kor}[theorem]{Corollary}

\newcommand{\cCS}{{\mathcal{CS}}}

\newcommand{\Pfaff}{{\tt Pfaff}}
\newcommand{\Map}{{\tt Map}}

\newcommand{\rk}{{\rangle}}
\newcommand{\lk}{{\langle}}

\newcommand{\Z}{\mathbb{Z}}

\newcommand{\diag}{{\tt diag}}
\newcommand{\proof}{{\it Proof.$\:\:\:\:$}}

\newcommand{\R}{\mathbb{R}}

\newcommand{\Q}{\mathbb{Q}}

\renewcommand{\det}{{\tt det}}
\renewcommand{\sinh}{\mathrm{sinh}}

\newcommand{\bL}{{\bf L}}

\newcommand{\C}{\mathbb{C}}
\renewcommand{\Pr}{{\tt Pr}}

\newcommand{\Tr}{{\tt Tr}}

\newcommand{\Line}{{\tt Line}}

\newcommand{\cH}{\mathcal{H}}

\newcommand{\cE}{\mathcal{E}}

\newcommand{\hol}{{\tt hol}}

\newcommand{\cW}{\mathcal{W}}

\newcommand{\cI}{\mathcal{I}}

\newcommand{\cU}{\mathcal{U}}
\newcommand{\Hom}{{\tt Hom}}

\newcommand{\End}{{\tt End}}

\newcommand{\im}{{\tt im}}

\newcommand{\ee}{{\tt e}}

\newcommand{\tr}{{\tt tr}}

\newcommand{\coker}{{\tt coker}}
\newcommand{\id}{{\tt id}}

\newcommand{\nat}{\mathbb{N}}

\newcommand{\cB}{\mathcal{B}}
\def\imath{{i}}

\def\hB{\hspace*{\fill}$\Box$ \newline\noindent}

\newcommand{\ind}{{\tt index}}

\newcommand{\cS}{\mathcal{S}}

\def\hB{\hspace*{\fill}$\Box$ \\[0.5cm]\noindent}

 \newcommand{\cG}{\mathcal{G}}

\newcommand{\cP}{\mathcal{P}}
 \newcommand{\cX}{\mathcal{X}}

\newcommand{\bW}{\mathbf{W}}
\newcommand{\pr}{{\tt pr}}

\newcommand{\ch}{{\mathbf{ch}}}

\newcommand{\bV}{\mathbf{V}}

\newcommand{\hHZ}{{\widehat{H\Z}}}

\begin{document}\maketitle

\begin{abstract}The present paper is a contribution to categorial index theory. Its main result is the calculation of the Pfaffian line bundle of a certain family of real Dirac operators as an object in the category of line bundles. Furthermore, it is  shown how string structures give rise to trivialisations of that Pfaffian.\end{abstract}

\tableofcontents
 
\section{Introduction}

\subsection{Topological, geometric and categorial aspects of index theory}

The present paper is a contribution to categorial index theory. Its main result is the calculation of the Pfaffian line bundle of a certain family of real Dirac operators as an object in the category of line bundles. Furthermore, we  show how string structures give rise to trivialisations of that Pfaffian.

Before we describe the results of the paper in greater detail in Subsection \ref{sss3} let us review the different levels  of 
index theory appearing in the title of the present Subsection.

The index of a family of elliptic differential operators $D$ parametrised by a space $B$ is a $K$-theory class $\ind(D)\in K^0(B)$. It is the homotopy class of an associated family of Fredhom operators defined by functional analytic techniques. If $D$ has a kernel bundle, then alternatively we can view the index as the formal difference $\ind(D)=[\ker(D)]-[\coker(D)]\in K^0(B)$ of vector bundles. Index theory provides the tools to calculate the $K$-theory class $\ind(D)\in K^0(B)$ or its cohomological invariants like its Chern character $\ch(\ind(D))\in H^{ev}(B;\Q)$ in terms of the symbol of $D$ \cite{MR0236950}.

If $D:=D(\cE)$ is the family of Dirac operators associated to a geometric family $\cE$ on a smooth manifold $B$, then its index can be refined to a geometric object. If the kernel bundle of $D$ exists, then it has an induced hermitean metric and a metric connection \cite[Ch. 9]{MR1215720}. This geometric information is encoded in the differential $K$-theory index class $\widehat{\ind}(D)\in \hat K^0(B)$ which can be defined in general without any assumption on the existence of a kernel bundle \cite{bunke-20071}. The differential index theorem calculates
the class $\widehat{\ind}(D)\in \hat K(B)$ or its Chern character
$\hat \ch(\widehat{\ind}(D))\in \widehat{H\Q}^{ev}(B)$ in differential rational cohomology in differential geometric terms, see  \cite{freed-2009} and \cite{bunke-20071}.

Let $E$ be the the total space of  the underlying proper submersion $\pi:E\to B$ of the geometric family $\cE$. If $W\to E$ is a complex vector bundle equipped with a hermitean metric $h^W$ and a metric connection $\nabla^W$, then we can form the twisted Dirac operator $D(\cE\otimes \bW)$, where by the bold face letter $\bW:=(W,h^W,\nabla^W)$ we denote the geometric bundle. In a categorial refinement of index theory one would consider the index of $D(\cE\otimes \bW)$ as an object $\widetilde{\ind}(D(\cE\otimes \bW))$ in a certain category $\tilde K^0(B)$ and study the functor $\bW\mapsto \widetilde{\ind}(D(\cE\otimes \bW))$
from the category of geometric vector bundles on $E$.  At the moment
such a theory is only partially understood, and the present paper discusses a particular aspect of that idea. For related developments in the context of algebraic geometry see 
\cite{MR902592} and  \cite{MR1085257}.

The first Chern class $c_1(\ind(D))\in H^2(B;\Z)$ classifies the topological type of the determinant line bundle of $D$. If $D=D(\cE)$ is the Dirac operator  associated to a geometric family $\cE$ on $B$, then $\det(D)$ comes with a Quillen metric $h^{\det(D)}$ and the Bismut-Freed connection
$h^{\det(D)}$, see \cite[Ch. 9 and  10]{MR1215720} for details.
The isomorphism class of the  geometric line bundle
$\det(D)$ over $B$ is classified by the first differential  Chern class
$\hat c_1(\det(D))\in \widehat{H\Z}^2(B)$ \cite{MR827262}. It can be derived from the differential $K$-theory class $\widehat{\ind}(D)$    by the identity 
$\hat c_1(\det(D))=\hat c_1(\widehat{\ind}(D))$, see \cite{bunke-chern}, \cite{MR2191484}. The calculation of integral Chern classes like
$c_1(\ind(D))$ and its differential refinements $\hat c_1(\widehat{\ind}(D))$  is the contents of integral index theory, see \cite{madsen} for first steps in the topological case.

The characterisation of the  determinant line bundle $\det(D)$ as an object in the category $\Line(B)$ of geometric line bundles over $B$ is an aspect of categorial index theory. In the presence of a real structure $J$
commuting with $D$, as observed in  \cite{MR915823},  one can define a natural square root of the inverse of the determinant line bundle\footnote{Note that our determinant line bundle, following the conventions in  \cite{MR1215720}, is the inverse of the determinant line bundle in  \cite{MR915823}. The Pfaffians coincide.}, the Pfaffian bundle $\Pfaff(D,J)$.

The main goal of the present paper is the calculation of the object
$\Pfaff(D,J)\in \Line(B)$ in a special situation which is motivated by applications in mathematical physics, notably string theory. 

\subsection{Description of the results}\label{sss3}

We consider a bundle $\pi:E\to B$ of compact two-dimensional manifolds, or alternatively, a proper submersion $\pi$ such that $\dim(E)-\dim(B)=2$.
Let further be given a 
fibrewise Riemannian metric $g^{T^v\pi}$ an a complement
$T^h\pi\subset TE$ of the vertical bundle $T^v\pi:=\ker(d\pi)$. In \cite{freed-2009}
the pair $(g^{T^v\pi},T^h\pi)$ is called  a Riemannian structure on $\pi$ since it gives rise to a Levi-Civita connection $\nabla^{T^v\pi}$, see  \cite[Ch. 9]{MR1215720}. We finally assume a  spin structure on $T^v\pi$ which allows to define the spinor bundle $S(T^v\pi)$.
In  \cite{MR2191484} we subsumed  this collection of data into the notion of a geometric family $\cE=(\pi:E\to B,g^{T^v\pi},T^h\pi,S(T^v\pi))$. 
By $D(\cE)$ we denote the Dirac operator associated to the geometric family $\cE$.

The spinor bundle associated to a two-dimensional vector bundle (like $T^v\pi$) with spin-structure is $\Z/2\Z$-graded and has a quaternionic  structure $J$ (see \cite[2.2.6]{MR2443109}),  a parallel, anti-linear, anti-selfadjoint, and  odd bundle endomorphism which commutes with the Clifford multiplication.

Let $\bV:=(V,h^V,\nabla^V)$ be a real geometric   vector bundle over $E$. Then we can form the Dirac bundle $S(\cE)\otimes_\R \bV$. Since we tensor over the real numbers, the quaternionic structure extends to the tensor product.  We let $$\cE\otimes_\R \bV:=(\pi:E\to B,g^{T^v\pi},T^h\pi,S(T^v\pi)\otimes_\R \bV)$$ denote the induced geometric family and use the symbol $J$ also to denote the extended  quaternionic structure.

The composition $JD(\cE\otimes_\R \bV)$ is an anti-linear, anti-selfadjoint, and even operator. By $(JD(\cE\otimes_\R \bV))^+$ we denote the component acting on sections of $S(\cE\otimes_\R \bV)^+$.
The relative Pfaffian line bundle 
$$\Pfaff(\cE\otimes_\R \bV,J,rel):=\Pfaff((JD(\cE\otimes_\R \bV))^+)\otimes \Pfaff((JD(\cE))^+)^{-n}\ ,\quad  n:=\dim(V)$$ is a complex line bundle with connection functorially
associated to this data described above. Its construction (see e.g. \cite{MR915823},  \cite{MR1186039}, \cite{MR2207325}, \cite{MR1971377}) will be recalled in Subsection \ref{pg111} below. The square of the relative Pfaffian line bundle is isomorphic to the inverse of the relative determinant line bundle
$\det(D(\cE\otimes_\R\bV))\otimes \det(D(\cE))^{-n}$ as a geometric line bundle. 

The Pfaffian line bundle plays an important role in two-dimensional quantum field theory
 \cite{MR2207325}, where the functional integral of the action over the fermions can be interpreted as a section of the  Pfaffian line bundle. In order to interpret the action as a complex valued function it is important to construct trivialisations of the Pfaffian line bundle.

The construction of the  line bundle $\Pfaff(\cE\otimes_\R\bV,J,rel)$ is based on the analytic properties of the family of Dirac operators $D(\cE\otimes_\R\bV)$. So it is not obvious how additional topological  and differential geometric structures can lead to a trivialisation.

In Subsection \ref{ztgudqwdqwd} of present paper we give a functorial, differential geometric construction
of a geometric line bundle $\bL=(L,h^L,\nabla^L)$ under the assumption that $V$ has a spin structure which is fibrewise trivial. The isomorphism class of
$\bL$ is classified by the first differential Chern class
$\hat c_1(\bL)\in \hHZ^2(B)$. The spin-characteristic class
$\frac{p_1}{2}(V)\in H^4(E;\Z)$ also has a differential refinement   $\frac{\hat p_1}{2}(\bV)\in \hHZ^4(E)$. 
We refer to \cite{MR827262} for a first construction of differential integral cohomology groups (there called groups of differential characters) and of the differential refinements of characteristic classes. We will give an alternative and well-adapted to the present purpose description of these classes in Subsection \ref{ll111}. Essentially by construction we have 
$$\hat c_1(\bL)=\int_{E/B} \frac{\hat p_1}{2}(\bV)\in \hHZ^2(B)\ ,$$
see Lemma \ref{juduqwdqwdqd}.
We have
\begin{theorem}[Theorem \ref{mmm1}]\label{hjjhdasd}
If $V$ has a spin structure which is fibrewise trivial, then
there is a functorial isomorphism of geometric line bundles
$$\Pfaff(\cE\otimes_\R\bV,J,rel)\cong \bL\ .$$
\end{theorem}
The adjective "functorial" here and above refers to base change along smooth maps $B^\prime\to B$.
As a consequence we get
\begin{kor}\label{dwqdwqdiui}
If $V$ has a spin structure which is fibrewise trivial, then
\begin{equation}\label{uidqwdqwd}
\hat c_1(\Pfaff(\cE\otimes_\R\bV,J,rel))= \int_{E/B} \frac{\hat p_1}{2}(\bV)\ .
\end{equation}
\end{kor}
The underlying equality $$c_1(\Pfaff(\cE\otimes_\R\bV,J))= \frac{\hat p_1}{2}(V)\ ,$$
which has first been derived in \cite[Prop. 5.4]{MR915823},
can be considered as an example of an integral index theorem, as alluded to in \cite{bunke-chern} (see  \cite{madsen} for a more general example), and Equation (\ref{uidqwdqwd})
can be considered as its differential refinement. 

Our second result concerns trivialisations of the relative Pfaffian line bundle $\Pfaff(\cE\otimes_\R\bV,J,rel)$. We assume that $\dim(V)=n\ge 3$. We define the homotopy type
$BString(n)$ as the homotopy fibre of  $\frac{p_1}{2}:BSpin(n)\to K(\Z,4)$.
 A topological string structure for a spin bundle $V$ is a homotopy class of lifts
$$\xymatrix{&BString(n)\ar[d]\\B\ar@{.>}[ur]\ar[r]^V&BSpin(n)}\ .$$

In the present paper we use an equivalent, more geometric notion of a string structure as a trivialisation of the Chern-Simons gerbe. This
and the  notion of a geometric string structure $str$ of a geometric spin bundle $\bV$ has  been discussed in \cite{waldorf}.
We will explain details at length in Subsection \ref{str111}. In particular, to a geometric string structure $str$ we have an associated
$3$-form $H_{str}\in \Omega^3(E)$, see (\ref{zdqwdqwddddqwdq}).
\begin{theorem}[Theorem \ref{zuadddqwdqd}]\label{tzu}
A  geometric string structure $str$ on $\bV$
gives rise to a functorial unit-norm section $s_{str}\in C^\infty(B,L)$ with
$$\nabla^L \log s_{str}=2\pi i \int_{E/B}H_{str}\ ,$$
\end{theorem}

If we combine Theorem \ref{hjjhdasd} with Theorem \ref{tzu} we get the following consequence.
\begin{kor}\label{uifwefwef}
A   geometric string structure $str$ on $\bV$
gives rise to a functorial unit-norm section  $s_{str}\in C^\infty(B,\Pfaff(\cE\otimes_\R\bV,J,rel)$ such that
$$\nabla^{\Pfaff(\cE\otimes_\R\bV,J,rel)}\log s_{str}= 2\pi i \int_{E/B}H_{str}\ .$$
\end{kor}

This corollary was the original motivation of the present paper. It answers a question
by Stephan Stolz. 

 Let us finally explain a typical example, which is the application looked for by physicists. We consider a compact surface $\Sigma$ 
with a Riemannian metric $g^{T\Sigma}$ and a spin structure. Furthermore, we consider a Riemannian manifold $X$ and a smooth map   $B\to \Map(\Sigma,X)$. Technically, this is a smooth map $$E:=\Sigma\times B\stackrel{f}{\to} X\ .$$  We let $\pi:E\to B$ be the projection,
$T^h\pi:=TB\subset T(\Sigma\times B)$ be the canonical horizontal subspace, and
$g^{T^v\pi}$  and the spin structure on $T^v\pi$ be induced by $g^{T\Sigma}$ and the spin structure on $\Sigma$, respectively. In this way we get a geometric family $\cE$. The real vector bundle
is obtained by $\bV:=f^*\mathbf{TX}$, where the geometric bundle $\mathbf{TX}=(TX,g^{TX},\nabla^{TX})$  is given by the Riemannian geometry of $X$, in particular, $\nabla^{TX}$ is the Levi-Civita connection. 

We assume that $X$ has  a string structure. 
In this situation by Corollary \ref{dwqdwqdiui} the isomorphism class of the relative Pfaffian line bundle can be calculated by transgressing the differential Pontrjagin class of $\mathbf{TX}$.
$$\hat c_1(
\Pfaff(\cE\otimes_\R\bV,J,rel))=\int_{\Sigma\times B/B} f^* (\frac{\hat p_1}{2}(\mathbf{TX}))\ .$$ Note that the string structure
on $X$ ensures that the spin structure of $f^*TX$ is fibre-wise trivial.
Furthermore, a refinement of the string structure of $X$ to a geometric string  structure $str$ gives rise to a trivialisation $s_{str}$ of $\Pfaff(\cE\otimes_\R\bV,J,rel)$ by Corollary \ref{uifwefwef}.

\medskip

{\em Acknowledgement: I thank Stephan Stolz for suggesting this  problem and Dan Freed for valuable hints.}

\section{The bundle $\bL$}\label{ll}

\subsection{Trivializations of $Spin$-structures and $\eta^3$-forms}\label{lllewuewe}

A $\Z/2\Z$-graded complex geometric vector bundle $\bW=(W,h^W,\nabla^W)$ over a manifold $M$ gives rise to a geometric family $\cW$ with zero-dimensional fibres, see \cite[2.2.2.1]{MR2191484}.  An invertible odd bundle endomorphism $\bar Q\in \End(W)^{odd}$ can be considered as a  taming $\cW_{t}$ in the sense of \cite[Def 2.2.4]{MR2191484}. To the tamed geometric family $\cW_t$ by   \cite[(2.26)]{MR966608} we can associate an eta form $\eta(\cW_t)\in \Omega^{odd}(M)$ such that $d\eta(\cW_t)=\ch(\nabla^{W})$.  We refer to  
(\ref{zudzhqwdwqdqwdqdwqdwqdw}) for our normalisations.  In the present subsection we study special properties
of the eta form in the case that $\bW$ comes from a real spin vector bundle $\bV$ and $\bar Q$ is induced by a spin trivialisation $Q$ of $V$. The main result is Lemma \ref{udqidwddwqdwqdwq}.

Let $\bV=(V,h^V,\nabla^V)$ be a real $n\ge 3$-dimensional geometric vector bundle over some manifold $M$. Then we form the $\Z/2\Z$-graded complex geometric bundle $\bW=(W,h^W,\nabla^W)$ such that
$\bW^+\cong \bV\otimes_\R\C$ and $\bW^-$ is the $n$-dimensional trivial geometric bundle over $M$.  A spin structure on $V$ is given by a $Spin(n)$-reduction   $Spin(V)\to O(V)$ of
the orthonormal frame bundle $O(V)$ of $V$. A trivialisation of the spin bundle $V$ is a trivialisation $Q:Spin(V)\stackrel{\sim}{\to} M\times Spin(n)$. The trivialisation of the spin bundle $V$ naturally gives rise to an unitary  vector bundle isomorphism
$Q^+:W^+\to W^-$.  We define the unitary odd involution
$$\bar Q:=\left(\begin{array}{cc}
0&(Q^+)^*\\Q^+&0
\end{array}\right)\in \End(W)^{odd}\ .$$
Let $\cW$ be the geometric family given by the $\Z/2\Z$-graded complex geometric vector bundle $\bW$.
 The local index form $\Omega(\cW)$ \cite[Def. 2.2.8]{MR2191484} in this case is the Chern-Weil representative $\ch(\nabla^W)\in \Omega^{ev}(M)$ of the Chern character of $W$ for the connection $\nabla^W$:
$$\Omega(\cW)=\ch(\nabla^W)=\ch(\nabla^{V\otimes_\R\C})-n\ .$$
We are in particular interested in the degree-$4$-component. We have
$c_1(\nabla^{V\otimes_\R\C})=0$ and therefore $$\ch_2(\nabla^{V\otimes_\R\C})=c_2(\nabla^{V\otimes_\R\C})=-p_1(\nabla^V)\ ,$$
where $c_i(\nabla^{V\otimes_\R\C})\in \Omega^{2i}(M)$ and $p_i(\nabla^{V})\in \Omega^{4i}$ again denote the Chern-Weil representatives of the corresponding characteristic classes.
 
The endomorphism $\bar Q$ gives rise to a taming $\cW_{t_Q}$  and an associated $\eta$-form. Its definition employs the rescaled
  super connection 
\begin{equation}\label{uzdzqwuidhziqwdq}
A_t:=t^{\frac{1}{2}} \bar Q+\nabla^W\ .
\end{equation} 
\begin{ddd}
The degree-$2k-1$ component of the eta-form of $\cW_{t_Q}$ is defined by
\begin{equation}\label{zudzhqwdwqdqwdqdwqdwqdw}
\eta^{2k-1}(\cW_{t_Q})=\frac{1}{(2\pi i)^k}\int_0^\infty \tr[ \partial_tA_t \exp(-A_t^2)]_{2k-1}\ .
\end{equation}\end{ddd}
Here $\eta^{2k-1}(\cW_{t_Q})\in \Omega^{2k-1}(M)$ denotes the component of the eta-form in degree $2k-1$, and $[\omega]_{2k-1}$ takes the degree-$2k-1$-component of the inhomogeneous form $\omega$.
Furthermore, if $W$ is a $\Z/2\Z$-graded vector bundle, then  $\tr:\End(W)\to \C$ denotes the super trace. 
The definition \cite[Def 2.2.16]{MR2191484} of the eta form of a tamed geometric family  is based on different family of super connections. It is related to (\ref{uzdzqwuidhziqwdq}) by a transformation  in the scaling parameter $t$.
Hence, it  gives the same eta form as in the present paper.

In the present paper we are in particular interested in
$$\eta^3(\cW_{t_Q})\in \Omega^3(M)\ .$$
It satisfies \begin{equation}\label{uidqwdqwd1}
d\eta^3(\cW_{t_Q})=\Omega^4(\cW)=-p_1(\nabla^V)\ .
\end{equation}

Let us calculate the eta form explicitly in order to see that it is nothing else then a classical Chern-Simons form.
Note that $W^-$ is trivialised so that we can identify sections of $\End(W^-)$ with matrix-valued functions.
 We define the matrix-valued one-form $B^+$ such that
$Q^+\nabla^{W^+}(Q^+)^*=d+B^+$, and  
the matrix-valued curvature two-form
$R^+:=dB^++B^{+2}$.
\begin{ddd}
The Chern-Simons form\footnote{The standard normalisation of the  Chern-Simons form is $$CS_{st}(\nabla^{W^+})=\tr R^+B^+-\frac{1}{3}\tr B^{+3}$$ such that $$d CS_{st}(\nabla^{W^+})=\tr (R^{\nabla^{W^+}})^2=2 (2\pi i)^2 c_2(\nabla^{W^+})\ .$$  The factor $\frac{1}{4}\frac{1}{(2\pi i)^2}$ is introduced for convenience.}   of the connection $\nabla^{W^+}$ (in the trivialisation
given by $Q^+$\footnote{One should better write $CS(\nabla^{W^+},Q^+)$, but we refrain from doing so in order to shorten the notation.}) is defined by 
$$ CS(\nabla^{W^+}):=\frac{1}{4}\frac{1}{(2\pi i)^2}\left(\tr R^+B^+-\frac{1}{3} \tr B^{+3}\right)\ .$$
  \end{ddd}
Its differential is  half of the second Chern form: 
$$d   CS(\nabla^{W^+})=  \frac{1}{2}c_2(\nabla^{W^+})\ .$$
\begin{lem}\label{uqifqffff77wefefwf}
We have
$$\eta^3(\cW_{t_Q})=2CS(\nabla^{W^+})\ .$$
\end{lem}
\proof
Define
$U:=\diag(Q^+,1)$. Then
$$B_t:=UA_tU^*=d+B+t^{\frac{1}{2}} K\ , \quad K:=\left(\begin{array}{cc}
0&1\\1&0
\end{array}\right)\\ ,\quad B:=\left(\begin{array}{cc}
B^+&0\\0&0
\end{array}\right) $$
is a rescaled superconnection on the trivial bundle $W^-\oplus W^-$ isomorphic to $A_t$.
Using the notation
$$ R:=\left(\begin{array}{cc}
R^+&0\\0&0
\end{array}\right)\ ,\quad C:=\left(\begin{array}{cc}
0&B^+\\-B^+&0
\end{array}\right)$$ we get 
$$B_t^2=dB+B^2+t^{\frac{1}{2}}\{B,K\}+tK^2=R+t^{\frac{1}{2}}C+t\ .$$
Furthermore
$$\partial_t B_t=\frac{1}{2t^{\frac{1}{2}}}K\ .$$ 
We now calculate the $3$-form component:
\begin{eqnarray*}
[\partial_t B_t\ee^{-B_t^2}]_3&=&\frac{\ee^{-t}}{2t^{\frac{1}{2}}}K\left(\frac{t^{\frac{1}{2}}}{2} (RC+CR)-\frac{t^{\frac{3}{2}}}{6} C^3\right)\\
&=&\frac{\ee^{-t}}{4}K\left(\begin{array}{cc}
0&R^+B^+\\-B^+R^+&0
\end{array}\right)-\frac{t\ee^{-t}}{12}K\left(\begin{array}{cc}
0&-B^3\\B^3&0
\end{array}\right)\\
&=&\frac{\ee^{-t}}{4}\left(\begin{array}{cc}
R^+B^+&0\\0&-B^+R^+
\end{array}\right)-\frac{t\ee^{-t}}{12}\left(\begin{array}{cc}
B^{+3}&0\\0&-B^{+3}
\end{array}\right)
\end{eqnarray*}
 We get
$$\tr [\partial_t B_t\ee^{-B_t^2}]_3=\frac{\ee^{-t}}{2} \tr R^+B^+
-\frac{t\ee^{-t}}{6} \tr B^{+3}\ .$$  
Using $$\int_0^\infty \ee^{-t}dt=1\ ,\quad \int_0^\infty t \ee^{-t}dt=1$$
we get
$$\eta^3(\cW_{t_Q})=\frac{1}{(2\pi i)^2} \left(\frac{1}{2}\tr R^+B^+-\frac{1}{6} \tr B^{+3}\right)=2CS(\nabla^{W^+})\ .$$
\hB

We now consider a second spin trivialisation $Q^\prime$ of the spin bundle $V$. An element $x\in H^3(M;\R)$ is called even, if it belongs to $2\im(H^3(M;\Z)\to H^3(M;\R))$.
\begin{lem}\label{udqidwddwqdwqdwq}
The difference
$\eta^3(\cW_{t_Q})-\eta^3(\cW_{t_{Q^\prime}})$ is closed and represents an even class in $H^3(M;\R)$.
\end{lem}
\proof
The difference is closed by (\ref{uidqwdqwd1}). That it represents an even class relies on the fact that we consider pairs of spin trivialisations as opposed to
just trivialisations. In order to  see this we invoke some index theory.
Alternatively one could use Lemma \ref{uqifqffff77wefefwf}
and show the corresponding property for the Chern-Simons form, a fact which is surely known to specialists, but for which we do not know a good reference.

We form a geometric family $\cI$ over $M$. The underlying fibre bundle of $\cI$  is
$p:[0,1]\times M\to M$ with the standard metric $p^{T^vp}$ and horizontal distribution $T^hp=TM\subset T([0,1]\times M)$. The trivial bundle $T^vp$ has an induced spin structure.
Note that $\partial (\cI\otimes p^* \bW)\cong \cW\oplus \cW^{op}$.
The tamings $\cW^{op}_{t_Q}$ and $\cW_{t_{Q^\prime}}$ induce a boundary taming
$(\cI\otimes \bW)_{bt}$. Since $\cI$ is an one-dimensional boundary tamed geometric family its index $\ind((\cI\otimes \bW)_{bt})$ is an element of $K^{-1}(M)$. We let $\ch^{odd}:K^{-1}(M)\to H^{odd}(M;\Q)$ be the odd Chern character.
By the index theorem for boundary tamed families \cite[Thm. 2.2.18]{MR2191484} we have the following equality in de Rham cohomology:
$$\ch^{odd}_3(\ind((\cI\otimes \bW)_{bt}))=[\Omega^3(\cI\otimes \bW)+\eta^3(\cW_{t_{Q^\prime}})-\eta^3(\cW_{t_Q})]\ .$$
We now observe that $$\Omega^3(\cI\otimes \bW)=\int_{[0,1]\times M/M} p^* \ch_2(\nabla^W)=0\ .$$
It remains to show that 
$\ch^{odd}_3(\ind(\cI\otimes \bW)_{bt})$ is even.
 
We first consider the special case
where $V_{Spin(n)}=Spin(n)\times \R^n$ is the trivial bundle over $Spin(n)$ with the trivial  metric, connection and spin structure $Spin(V_{Spin(n)})=Spin(n)\times Spin(n)$. For $Q_{Spin(n)}$ we take the canonical trivialisation, and for $Q_{Spin(n)}^\prime$ we take the universal trivialisation $Q_{Spin(n)}^\prime (h,g)=(h,hg)$.  This defines by the construction above a boundary tamed geometric family $(\cI_{Spin(n)}\otimes \bW_{Spin(n)})_{bt}$ over $Spin(n)$.
We first claim that $\ch^{odd}_3(\ind(\cI_{Spin(n)}\otimes \bW_{Spin(n)})_{bt})\in H^3(Spin(n);\R)$ is even. In order to see this we show that $(\cI_{Spin(n)}\otimes \bW_{Spin(n)})_{bt}$ is a pull-back of a boundary tamed family via the two-fold covering $s:Spin(n)\to SO(n)$. Indeed, we
can consider the trivial bundle 
$V_{SO(n)}\to SO(n)$ and $W_{SO(n)}:= V_{SO(n)}\oplus  V_{SO(n)}^{op}$. 
We set
$$\bar Q_{SO(n)}:=\left(\begin{array}{cc}
0&(Q_{SO(n)}^+)^*\\Q_{SO(n)}^+&0
\end{array}\right)\in \End(W_{SO(n)})\ ,$$
where $\bar Q_{SO(n)}(g,x)=(g,gx)$, $(g,x)\in SO(n)\times \R^n=V_{SO(n)}$.
As above this gives a boundary tamed geometric family
$(\cI_{SO(n)}\otimes \bW_{SO(n)})_{bt}$, and we have $$(\cI_{Spin(n)}\otimes \bW_{Spin(n)})_{bt}\cong s^*( \cI_{SO(n)}\otimes \bW_{SO(n)})_{bt}\ .$$
Note that in general $\ch^{odd}_3(x)\in H^3(X;\R)$ is integral for every $x\in K^{-1}(X)$. 
Therefore $\ch^{odd}_{3}(\ind(( \cI_{SO(n)}\otimes \bW_{SO(n)})_{bt}))$ is integral. 
Since $s^*:H^3(SO(n);\R)\to H^3(Spin(n);\R)$ maps integral to even elements, we conclude that
$$\ch^{odd}_3( \ind((\cI_{Spin(n)}\otimes \bW_{Spin(n)})_{bt}))=s^* \ch^{odd}_{3}(\ind( \cI_{SO(n)}\otimes \bW_{SO(n)})_{bt})$$ is even.

We now come back our original case. The transition
$Q^\prime\circ Q^{-1}:M\times Spin(n)\to M\times Spin(n)$  determines a map
$T:M\to Spin(n)$ such that
$Q^\prime\circ Q^{-1}(m,g)=(m,T(g)m)$.
We observe that $T^* (\cI_{Spin(n)}\otimes \bW_{Spin(n)})_{bt}$ differs, up to isomorphism induced by $Q$, from
$(\cI\otimes \bW)_{bt}$ only by the choice of the connection on the twisting bundle.
But the index as a homotopy invariant is independent of the connection so that
$$\ch^{odd}_3(\ind((\cI\otimes \bW)_{bt}))=T^*\ch^{odd}_3(\ind((\cI_{Spin(n)}\otimes \bW_{Spin(n)})_{bt}))$$
is an even class.
\hB

\subsection{The Cheeger-Simons character $\frac{\hat p_1}{2}(\bV)$
} \label{ll111}

A generalised cohomology theory has differential extensions. We refer to
\cite{bunke-2009} for a description of the axioms of a differential extension of a generalised cohomology theory and to \cite{MR2192936} for a general construction. The differential extension of complex $K$-theory already appeared in the introduction to the present paper.  But in detail we will only need the differential extension of integral cohomology  $(\hHZ^*,R,I,a)$, where
$$R:\hHZ^*(M)\to \Omega^*(M)\ , \quad  \hHZ^*(M)\to H\Z^*(M)\ ,\quad a:\Omega^{*-1}(M)\to \hHZ^*(M)$$
denote the structure maps. 
By \cite{MR2365651} or \cite{bunke-2009}
the differential extension  of integral cohomology is unique and multiplicative.

In the original model of  \cite{MR827262} the $k$'th
differential integral cohomology  group $\hHZ^k(M)$
is the group of differential characters of degree $k$.  Let $Z_{k-1}(M)$ denote the group of smooth $k-1$-dimensional cycles in $M$.
A differential character of degree $k$
is a homomorphism $\phi:Z_{k-1}(M)\to U(1)$ such that there exists a smooth form
$R(\phi)\in \Omega^k(M)$ so that
$\phi(\partial c)=\exp\left(2\pi i \int_{c} R(\phi)\right)$
for all smooth chains $c\in C_k(M)$. It suffices to verify this condition for smooth simplices.

A real spin vector bundle $V\to M$ has a characteristic class
$\frac{p_1}{2}(V)\in H^4(M;\Z)$. For a geometric spin bundle
$\bV$ this characteristic  class has a differential refinement 
$\frac{\hat p_1}{2}(\bV)\in \hHZ^4(M)$ first defined  using Chern-Weil theory in \cite{MR827262}. 
 We now describe this class as a differential character $\frac{\hat p_1}{2}(\bV):Z_3(M)\to U(1)$.
If $z\in Z_3(M)$, then we can choose a neighbourhood $U$ of the trace $|z|$ of $z$ which is homotopy equivalent to a three-dimensional $CW$-complex. Furthermore, since
$BSpin(n)$ is $3$-connected, we can choose a trivialisation
$Q$ of the spin bundle $V_U$. \begin{ddd}We define
\begin{equation}\label{uwidqwuidhqwdqwdqdqd}
\frac{\hat p_1}{2}(\bV)(z):=\exp\left(-\pi i \int_z \eta^3(\cW_{t_Q})\right)\ .
\end{equation}
\end{ddd}
We must show that this does not depend on the choice of 
the trivialisation. Indeed, if  $Q^\prime$ is a second trivialisation,
then by Lemma  \ref{udqidwddwqdwqdwq}.
$$ \int_z \eta^3(\cW_{t_{Q^\prime}}) - \int_z \eta^3(\cW_{t_Q})\in 2\Z\ .$$

If $c\in C_4(M)$ is a smooth $4$-simplex, then still we can choose a trivialisation of $Q$ of the spin bundle $V_U$ on some neighbourhood $U$ of $|\sigma|$. In this case we get by Stoke's theorem and (\ref{uidqwdqwd1})
$$\frac{\hat p_1}{2}(\bV)(\partial c)=\exp\left(-\pi i \int_c d\eta^3(\cW_{t_{Q }}\right)= \exp\left(\pi i \int_{c} p_1(\nabla^V)\right)\ .$$
Therefore $$R(\frac{\hat p_1}{2}(\bV))=\frac{1}{2} p_1(\nabla^V)$$
as it should be.

\begin{lem}
The differential cohomology class defined by the differential character above coincides with the class $\frac{\hat p_1}{2}(\bV)$ defined in \cite{MR827262}.
\end{lem}
\proof
For the moment, let us denote the differential cohomology class corresponding to the differential character described above by $\hat \phi(\bV)$.
It is easy to see that 
 $\hat \phi(\bV)$ is natural with respect to base change along smooth maps $M\to M^\prime$. Since $BSpin(n)$ is $3$-connected and $H^4(BSpin(n);\Z)$ is torsion-free 
there exists a smooth manifold $M^\prime$ with a real $n$-dimensional geometric spin bundle $\bV^\prime$ and a map $f:M\to M^\prime$ such that
$\bV\cong f^*\bV^\prime$, this isomorphism is covered by 
an isomorphism $Spin(V)\cong f^*Spin(V^\prime)$, and such that
$H^3(M^\prime;\R)=0$ and $H^4(M^\prime;\Z)$ is torsion free.
Under these conditions a class
$\hat x \in \hHZ^4(M^\prime)$ is completely determined by
its curvature $R(\hat x)\in \Omega^4(M^\prime)$. Since $R(\hat \phi(\bV))=\frac{1}{2}p_1(\nabla^{V^\prime})=R(\frac{\hat p_1}{2}(\bV^\prime))$
we have
$\hat \phi(\bV^\prime)=\frac{\hat p_1}{2}(\bV^\prime)$.
By naturality this implies 
$\hat \phi(\bV)=\frac{\hat p_1}{2}(\bV)$.
\hB


\subsection{The construction of $\bL$}\label{ztgudqwdqwd}

Let $\pi:E\to B$ be a bundle of compact oriented  two-dimensional manifolds.
In the present subsection we construct a geometric line bundle
$\bL=(L,h^L,\nabla^L)$ over $B$ functorially associated to a geometric spin bundle $\bV$ over $E$ which is trivial as a spin bundle over the fibres of $\pi$.
More precisely, by functoriality we mean that for a smooth map
$f:B^\prime \to B$ and induced cartesian diagram
\begin{equation}\label{uqidqwdqwdwqdqd545454545}
\xymatrix{E^\prime\ar[r]^F\ar[d]^{\pi^\prime}&E\ar[d]^{\pi}\\B^\prime\ar[r]^f&B}
\end{equation} we have an associated isomorphism
$\phi_f:f^*\bL\stackrel{\sim}{\to} \bL^\prime$, where $\bL^\prime$ is the line bundle associated to $F^*\bV$, and for a second smooth map
$g:B^{\prime\prime}\to B^\prime$ we have the associativity relation
$\phi_g\circ g^*\phi_f=\phi_{f\circ g}$.

 We are going to describe the geometric bundle $\bL$ by describing its local sections. We first describe a set-valued presheaf
$I_{\dots}:B\supseteq U\mapsto I_U$ on $B$ such that
$I_U$ is non-empty if $U$ is contractible.
The elements $Q\in I_Q$ will index local sections $s_Q\in C^\infty(U,L)$.
After all, the maps $I_U\ni Q\mapsto s_Q\in 
C^\infty(U,L)$ for all $U\subseteq B$ together combine to  a map $I_{\dots}\to C^\infty(\dots,L)$  of presheaves from $ I_{\dots}$ to the sheaf of local
sections of $L$, but this can be said only after we know what $L$ is.

In order to complete the definition of  $(L,h^L)$ we must provide the transition
functions $$\frac{s_{Q^\prime}}{s_Q}=c(Q^\prime,Q)\in C^\infty(U,U(1))$$ for pairs
$Q,Q^\prime\in I_U$, see (\ref{uiwqddqwd44}). These functions must fit into a map
of presehaves $I_{\dots}\times I_{\dots}\to C^\infty(\dots,U(1))$
and satisfy a cocycle relation. In order to construct
the connection $\nabla^L$ of $L$ we will explicitly (see (\ref{uiwqddqwd444}))
    describe  the connection one-forms
$$\omega_Q:=\nabla^{L}\log s_Q:=s_Q^{-1}\nabla^Ls_Q\in \Omega^1(U)$$ 
and verify their compatibility with the transition functions.
 
We let $I_U$ be the set trivialisationsons $Q:Spin(V_{E_U})\stackrel{\sim}{\to} E_U\times Spin(n)$, where $E_U=\pi^{-1}(U)$.
Recall that that we assume that the spin structure $Spin(V)\to E$ is trivial on the fibres of $\pi:E\to B$. Hence, if  $U\subseteq B$ is  contractible, then  
$I_U\not=\emptyset$.

The spin trivialisation 
 $Q\in I_U$ of   $V_{E_U}\to E_U$ induces  a tamed geometric family
$\cW_{t_Q}$ over $E_U$ as described in  Subsection \ref{lllewuewe}.

We define the connection one-form 
\begin{equation}\label{uiwqddqwd444}
\omega_Q:=-\pi i\int_{E_U/U} \eta^3(\cW_{t_Q})\in \Omega^1(U)\ .
\end{equation}
 
Since it is imaginary it defines a metric connection.

It suffices to define the transition functions
for contractible open subsets $U\subseteq B$
provided that they are compatible with the restriction to contractible
$U^\prime\subseteq U$, see (\ref{ioeufoewffewfwefwf}).

Assume that $U$ is contractible
and consider 
 $Q,Q^\prime\in I_{U}$. Since the fibres of $\pi:E\to B$ are two-dimensional, the manifold $E_{U}$ is homotopy equivalent to an at most  two-dimensional $CW$-complex.
Since $Spin(n)$ is two-connected, there exists a trivialisation $H$ of Spin bundles $\pr^*V_{E_{ U}}$ connecting  $Q$ with
$Q^\prime$, where $\pr:[ 0,1]\times E_{U}\to E_{U}$ is the projection.

We define
\begin{equation}\label{uiwqddqwd44}
c(Q^\prime,Q):=\exp\left(-\pi i \int_{[0,1]\times E_{ U}/U} \eta^3(\pr^*\cW_{t_{H}})\right)\in C^\infty(U,U(1))\ .
\end{equation}
This is independent of the choice of $H$. Indeed, if $H^\prime$ is another choice, then
we can concatenate these two choices in order to get a $Spin$-trivialisation
$G:=H\sharp H^{op, \prime}$ of $\pr^*V_{E_{U}}$, where $\pr$ is now the projection $\pr:S^1\times  E_{U}\to E_{U}$.
We must show that
$$\int_{S^1\times E_{U}/ U} \eta^3(\pr^*\cW_{t_G})\in C^\infty(U,2\Z)\ .$$  
Let $u\in U$, $E_u:=\pi^{-1}(u)$ be the fibre over $u$, and let $q
:S^1\times E_u\to E$ be the restriction of $\pr$. Since $q$ factors over the two-dimensional manifold $E_u$ we conclude that $0=q^*:\hHZ^4(E)\to \hHZ^4(S^1\times E_u)$.
By the construction of the differential character $\frac{\hat p_1}{2}(\bV)$ in Subsection \ref{ll111} we have
$$\frac{\hat p_1}{2}(q^*\bV)=a(-\frac{1}{2}\eta^3(q^*\cW_{t_G}))=q^*a(-\frac{1}{2}\eta^3(\cW_{t_G}))=0\ ,$$
where $a:\Omega^3(\dots)\to \hHZ^4(\dots)$ denotes one of the structure maps of the differential extension $\hHZ^4$.
This implies that
$$\left(\int_{S^1\times E_{U}/ U} \eta^3(\pr^*\cW_{ t_G})\right)(u)\in 2\Z\ .$$
It immediately follows from the construction that
\begin{equation}\label{ioeufoewffewfwefwf}
c(Q^\prime,Q)_{|U^\prime}=c(Q^\prime_{|U^\prime},Q_{|U^\prime})
\end{equation}  for a contractible open subset $U^\prime\subseteq U$.

We now check the cocycle condition. Let $Q^{\prime\prime}\in I_{U }$ be a third trivialisation. Then we can concatenate  the path $H$ with a path $H^\prime$ from $Q^\prime$ to $Q^{\prime\prime}$ to a path 
$G:=H^\prime \sharp H  $ from $Q$ to $Q^{\prime\prime}$.
The cocycle condition now follows from
$$\int_{[0,1]\times E_U/U}\eta^3(\pr^*\cW_{t_H})+\int_{[0,1]\times E_U/U}\eta^3(\pr^*\cW_{t_{H^\prime}})=\int_{[0,1]\times E_U/U}\eta^3(\pr^*\cW_{t_G})\ .$$

Finally we check the compatibility with the connection one-forms.
We must show, that
$$\omega_{Q^\prime}-\omega_Q=c(Q^\prime,Q)^{-1}dc(Q^\prime,Q)\ .$$
Indeed, by Stoke's theorem it follows immediately from the definition (\ref{uiwqddqwd44})
and  $$\int_{[0,1]\times E_{ U}/  U}d\eta^3(\pr^*\cW_{   t_{H}} )=\int_{[0,1]\times E_{  U}/  U}\pr^* p_1(\nabla^V)=0$$ that
$$c(Q^\prime,Q)^{-1}dc(Q^\prime,Q)= -\pi i\left(\int_{ E_U/U} \eta^3(\cW_{t_{Q^\prime}})-\int_{ E_U/U}\eta^3(\cW_{t_Q})\right)\ .$$ 
This finishes the construction of the bundle $\bL$.

Given a cartesian diagram (\ref{uqidqwdqwdwqdqd545454545}) and an open subset $U\subseteq B$, by pulling back trivialisations we get a map
$F^*:I_U\to I_{f^{-1}(U)}$. The bundle map $\phi_f:f^*L\to L^\prime$ is characterised by the property that for $Q\in I_U$  it maps the section $f^*s_Q\in C^\infty(f^{-1}(U),f^*L)$ to the section 
 $s_{F^*Q}\in C^\infty(f^{-1}(U),L^\prime)$.
One easily checks, using the functoriality of $\eta$-forms, that this defines an isomorphism of geometric line bundles
which behaves as required under compositions.

If $\pi:E\to B$ is a proper submersion with an orientation of the vertical bundle $T^v\pi$, then there is an integration map
$$\int_{E/B}:\hHZ^{*}(E)\to \hHZ^{*-\dim(E)+\dim(B)}(B)\ .$$
It has a very convenient description in the model introduced in \cite{bks}, see also the literature cited therein.
In the present paper we will only need the existence of the integration map,
its compatibility with cartesian diagrams of the form (\ref{uqidqwdqwdwqdqd545454545}) in the sense that
$$\int_{E^\prime/B^\prime} F^*\hat x=f^* \int_{E/B}\hat x$$
for $\hat x\in \hHZ^*(E)$, 
 and the property, that for  $\alpha\in \Omega^{*-1}(E)$ we have
$$\int_{E/B} a(\alpha)=a(\int_{E/B} \alpha)\ .$$

We now come back to our surface bundle $\pi:E\to B$.
We calculate the first differential Chern class $\hat c_1(\bL)\in \hHZ^2(B)$ of the geometric line bundle $\bL$. This task is equivalent to the calculation of the holonomy of $\bL$ since
the corresponding differential character $\hat c_1(\bL):Z_1(B)\to U(1)$
associates to the smooth cycle $z\in Z_1(B)$ the holonomy $ \hol(\bL)(z)\in U(1)$ of $\bL$  along $z$.

\begin{lem}\label{juduqwdqwdqd}
We have
$$\hat c_1(\bL)=\int_{E/B} \frac{\hat p_1}{2}(\bV)\ .$$
\end{lem}
 \proof
Let $z\in Z_1(B)$ be a smooth one-cycle. 
Then we can find a neighbourhood $U$ of its trace $|z|$ which is homotopy equivalent to a one-dimensional $CW$-complex. Since  $BSpin(n)$ is $3$-connected and  $E_U:=\pi^{-1}(U)$
is homotopy equivalent to a three-dimensional $CW$-complex  there exists an element $Q\in I_U$.
But then
$$\frac{\hat p_1}{2}(\bV)_{|E_U}=a(-\frac{1}{2}  \eta^3(\cW_{t_Q}))\ .$$
Hence
$$ \int_{E_U/U}\frac{\hat p_1}{2}(\bV)_{|E_U}=a(-\frac{1}{2}\int_{E_U/U}  \eta^3(\cW_{t_Q}))$$
so that
$$\left(\int_{E_U/U}\frac{\hat p_1}{2}(\bV)\right)(z) =\left(\int_{E_U/U}\frac{\hat p_1}{2}(\bV)_{|E_U}\right)(z)=\exp\left(-\pi i \int_z \int_{E_U/U}  \eta^3(\cW_{t_Q})\right)\ .$$
The spin trivialisation $Q$ of $V_{E_U}$ gives rise to the section $s_Q$  
and the corresponding connection one-form $\omega_Q$ defined by (\ref{uiwqddqwd444}).
We have
$$\hol(\bL)(z)=\exp(\int_z \omega_Q)=\exp\left(-\pi i\int_z \int_{E_U/U} \eta^3(\cW_{t_Q})\right) \ .
$$
Hence
$$\hat c_1(\bL)(z)=\hol(\bL)(z)=\left(\int_{E_U/U}\frac{\hat p_1}{2}(\bV)\right)(z)\ .$$
\hB

\section{The Pfaffian}\label{pg}

\subsection{Pfaffian and determinant line bundles for families of Dirac operators}
\label{pg111}
In this subsection we recall the construction of the Pfaffian line bundle of a family of Dirac operator with a real structure.
Let $\cE$ be an even  geometric family over a base $B$
such that the underlying Clifford bundle has an odd anti-linear, anti-selfadjoint
automorphism $J$. In particular, $J$ is parallel and commutes with the Clifford multiplication.  In this situation we define a geometric line bundle
$\Pfaff(\cE,J)$ over $B$. It is functorial and comes with a canonical isomorphism
\begin{equation}\label{udidwqdqwdwqdqdwq}
\kappa:\Pfaff(\cE,J)^2\stackrel{\sim}{\to} \det(\cE)^{-1}\ ,
\end{equation}
where $ \det(\cE)$ is the determinant line bundle of $\cE$, see e.g. \cite[Ch 9 and 10]{MR1215720} and the recapitulation below. In order to construct $\Pfaff(\cE,J)$ we will first construct the underlying complex line bundle together with the isomorphism   (\ref{udidwqdqwdwqdqdwq}). We then define the geometry on $\Pfaff(\cE,J)$ in the unique way such that  (\ref{udidwqdqwdwqdqdwq}) becomes an isomorphism of geometric bundles. 

First we recall the construction  \cite{MR915823}, \cite{MR1186039}, \cite{MR2207325}. The geometric family
$\cE$ gives rise to a family of Dirac operators $D(\cE)$ which acts on the bundle of $\Z/2\Z$-graded Hilbert spaces $H(\cE)$. We consider the families of anti-linear and anti-selfadjoint operators $D^\pm:=(JD(\cE))^\pm$ which act on the subbundles $H(\cE)^\pm$. The compositions $\Delta^\pm:=D^\mp(\cE) D(\cE)^\pm$ are non-negative and selfadjoint. For $\lambda\in [0,\infty)$ we consider the open subset $$U_\lambda:=\{u\in B\:|\:\lambda^2\not\in \sigma(\Delta^+(u))\cup   \sigma(\Delta^-(u))\}\subseteq B\ ,$$  where $\sigma(\Delta^\pm(u))$ denotes the spectrum of the operator $\Delta^\pm(u)$ over the point $u\in B$. By $E_{\Delta^\pm}[0,\lambda]$ we denote spectral projection of $\Delta^\pm$ onto the interval $[0,\lambda]$.  This family of projections is smooth over $U_\lambda$. Its image $H^\pm_\lambda:=E_{\Delta^\pm}[0,\lambda]H(\cE)^\pm$ is a finite-dimensional smooth complex vector bundle.

For an $n$-dimensional complex vector bundle $W\to X$ we let $\det(W)\to X$ denote the maximal non-trivial alternating power $\det(W) :=\Lambda^{n}W$. 

On $U_\lambda$ we define the line bundle
$$\Pfaff(\cE,J)_\lambda:=\det(H^+_\lambda)\ .$$
If $\mu>\lambda$, then
over $U_\lambda\cap U_\mu$ we have an orthogonal decomposition
$H^+_\mu\cong H^+_\lambda\oplus  H^+_{\lambda,\mu}$ of smooth bundles of Hilbert spaces, where
$H_{\lambda,\mu}^+=E_{\Delta^+}(\lambda,\mu]H(\cE)^+$.
On $U_\lambda\cap U_\mu$ we get an isomorphism 
$$\Pfaff(\cE,J)_\mu\cong \Pfaff(\cE,J)_\lambda\otimes \det(H^+_{\lambda,\mu})\ .$$
Note that
$D^+_{|H_{\lambda,\mu}}$ is an anti-linear anti-symmetric isomorphism of $H^+_{\lambda,\mu}$. Therefore the  form
$$d_{\lambda,\mu}(\dots,\dots):=\langle D^+_{|H_{\lambda,\mu}} \dots,\dots\rangle\in (\Lambda^2 H^+_{\lambda,\mu})^{*}$$
is nowhere vanishing. Hence we get a nowhere vanishing section
$$\Pfaff(D^+_{|H_{\lambda,\mu}}):=d_{\lambda,\mu}^{-\frac{\dim(H^+_{\lambda,\mu})}{2}}\in  C^\infty(U_\lambda\cap U_\mu,\det(H^+_{\lambda,\mu}))\ .$$
We define an isomorphism
$$c_{\lambda,\mu}:\Pfaff(\cE,J)_\lambda\stackrel{\sim}{\to} \Pfaff(\cE,J)_\mu\ ,\quad c_{\lambda,\mu}(s):= s\otimes \Pfaff(D^+_{|H_{\lambda,\mu}}) \ .$$
For $\nu\ge \mu\ge \lambda$ these isomorphisms satisfy the cocycle condition
$$c_{\mu,\nu}\circ c_{\lambda,\mu}=c_{\lambda,\nu}\ .$$ 
We glue the collection of line bundles
$\left(\Pfaff(\cE,J)_\lambda\to U_\lambda\right)_{\lambda\ge 0}$
using the cocycle $(c_{\lambda,\mu})_{\lambda,\mu\ge 0}$ in order to get the underlying  complex line bundle of $\Pfaff(\cE,J)$.
For the construction of the metric and the connection we invoke a canonical isomorphism
$$\kappa:\Pfaff(\cE,J)^2\stackrel{\sim}{\to} \det(\cE)^{-1}\ .$$
To this end we recall the very similar construction of the determinant line bundle $\det(\cE)$.
Over $U_\lambda$ we define the line bundle 
$$\det(\cE)_\lambda:=\det(H_\lambda^+)^{-1}\otimes \det(H_\lambda^-)\ .$$
For $\mu >\lambda$ on $U_\lambda\cap U_\mu$ we have isomorphisms
$$\det(\cE)_\mu\cong \det(\cE)_\lambda\otimes \det(H_{\lambda,\mu}^+)^{-1}\otimes \det(H_{\lambda,\mu}^-)\ .$$
The operator
$$D^+(\cE)_{|H_{\lambda,\mu}^+} :H_{\lambda,\mu}^+\to H_{\lambda,\mu}^-$$
is an isomorphism and therefore gives a nowhere vanishing section
$$\det(D^+(\cE)_{|H_{\lambda,\mu}^+})\in C^\infty(U_\lambda\cap U_\mu,
\det(H_{\lambda,\mu}^+)^{-1}\otimes \det(H_{\lambda,\mu}^-))\ .$$
We define the cocycle 
$$f_{\lambda,\mu}: \det(\cE)_\lambda\to \det(\cE)_\mu\ ,\quad f_{\lambda,\mu}(s):=s\otimes \det(D^+(\cE)_{|H_{\lambda,\mu}^+} )\ .$$
The underlying complex vector bundle of $\det(\cE)$ is obtained by glueing the family of bundles
$\left(\det(\cE)_\lambda\to U_\lambda\right)_{\lambda\ge 0}$ using the cocycle $(f_{\lambda,\mu})_{\lambda,\mu\ge 0}$.

We now define the isomorphism
$\kappa:\Pfaff(\cE,J)^2\stackrel{\sim}{\to}\det(\cE)^{-1}$ by defining a collection $(\kappa_\lambda)_{\lambda\ge 0}$ of isomorphisms
$$\xymatrix{\det(H_\lambda^+)\otimes \det(H_\lambda^+)\ar@{=}[d]\ar[r]^{\kappa_\lambda}_\cong&\det(H_{\lambda,\mu}^+)\otimes \det(H_{\lambda,\mu}^-)^{-1}\ar@{=}[d]\\\Pfaff(\cE,J)^2_\lambda\ar[r]_\cong& \det(\cE)_\lambda}$$
Indeed, the part $J^+$ of the anti-linear  isomorphism $J$ induces an isomorphism
$J^+_\lambda:H_\lambda^+\stackrel{\sim}{\to} \bar H_\lambda^-$, and therefore
an isomorphism $$\det(J_\lambda^+):\det(H_\lambda^+)\to \det(H_{\lambda}^-)^{-1}\ .$$
We set
$$\kappa_\lambda:=1\otimes \det(J_\lambda^+)\ .$$
It is easy to check that the collection $(\kappa_\lambda)_{\lambda\ge 0}$ is compatible with the cocycles and therefore defines an isomorphism $\kappa$ as required.
As explained above the geometry of $\Pfaff(\cE,J)$, i.e. the metric and the connection, is now defined in the unique way such that $\kappa$ becomes an isomorphism of geometric line bundles. This finishes the description of the Pfaffian bundle $\Pfaff(\cE,J)$ as a geometric line  bundle\footnote{We again remind the reader that our definition is the inverse of the Pfaffian 
in \cite{MR915823}}.

The family of Dirac operators $D(\cE)$ is invertible on $U_0$. Moreover, since
$H_0^\pm=0$ we have a canonical isomorphism
$$\det(H_0^\pm)\cong U_0\times \C$$ and therefore isomorphisms
$$\Pfaff(\cE,J)_0\cong U_0\times \C\ ,\quad \det(\cE)_0\cong U_0\times \C\ .$$
We let $$s^{1/2}_{can}\in C^\infty(U_0,\Pfaff(\cE,J))\ ,\quad s_{can}\in C^\infty(U_0,\det(\cE))$$
denote the corresponding sections.
Then 
$$\kappa(s^{1/2}_{can}\otimes s^{1/2}_{can})=s_{can}^{-1}\ .$$ 
By \cite[Ch 9]{MR1215720}, 
\begin{equation}\label{uqwidqwdwqdwd666}
\|s_{can}\|^2=\det(\Delta^+)\ ,\quad \|s^{1/2}_{can}\|^2=\sqrt{\det(\Delta^+)}^{-1}\ .
\end{equation}
We let $s_{can}^0$ and $s_{can}^{0,1/2}$ denote the normalised sections
$$s_{can}^0:= \frac{s_{can}}{\|s_{can}\|}\ ,\quad s_{can}^{1/2,0}:= \frac{s^{1/2}_{can}}{\|s^{1/2}_{can}\|}\ .$$
The connection one-form of $\det(\cE)$ is determined by 
$$ \nabla^{\det(\cE)} \log s^0_{can}=2\pi i \eta^1(\cE_{U_0,t})\ ,$$
where $\cE_{U_0,t}$ is the canonical taming of $\cE_{U_0}$ (by the zero perturbation).
Therefore the connection one-form for the Pfaffian bundle is given by 
\begin{equation}\label{udduqwdqwdqwd777}
 \nabla^{\Pfaff(\cE,J)} \log s^{1/2,0}_{can}
=- \pi i \eta^1(\cE_{U_0,t})\ .
\end{equation}
 Finally we calculate the curvature (see \cite[Thm 10.35]{MR1215720},  \cite[Thm 3.1]{MR915823}):
\begin{equation}\label{udquwdiqwudqwddq}
\quad  R(\det(\cE))= 2\pi i \Omega^2(\cE)\ ,\quad R(\Pfaff(\cE,J))=-\pi i \Omega^2(\cE)\ .
\end{equation}

\subsection{Theory for generalised Dirac operators}

The Dirac operator $D(\cE)$ of a geometric family is compatible, i.e. associated to a bundle of Clifford modules. This fact is important if one wants to use the standard local index theory calculations \cite{MR1215720} which e.g. give (\ref{udquwdiqwudqwddq}).  On the other hand, the construction of the Pfaffian and determinant line bundle works equally well for generalised Dirac operators, i.e. zero-order perturbations of compatible Dirac operators.
In the present paper we will consider generalised Dirac operators which arise as follows.
Recall that we consider a two-dimensional geometric family $\cE$ with underlying surface bundle $\pi:E\to B$ whose
Dirac bundle is  the spinor bundle $S(T^v\pi)$, and a real geometric vector bundle
$\bV$ over $E$. Our final goal is  study the Pfaffian of the family of Dirac operators associated to the geometric family $\cE\otimes \bW$, where $\bW:=\bW^+\oplus \bW^-$ with $\bW^+:=\bV\otimes_\R\C$ and $\bW^-$ is a trivial bundle of dimension $\dim_\R V$. The odd anti-involution $J$ is induced by the quaternionic structure of $S(T^v\pi)$ and the real structure of $W$.

If $\bar Q\in \End(W)$ is an odd,  selfadjoint endomorphism which commutes with the real structure of $W$, then we can form the family of generalised Dirac operator $$D(\cE\otimes \bW,\bar Q):=D(\cE\otimes \bW)+1\otimes \bar Q\ .$$
Its Pfaffian and determinant line bundles will be denoted by 
$$\Pfaff(\cE\otimes \bW,\bar Q,J)\ ,\quad \det(\cE\otimes \bW,\bar Q)\ .$$ 

Let us now discuss the contribution of the additional term $1\otimes \bar Q$ to the local index
calculation. For simplicity we just write $\bar Q$ instead of $1\otimes \bar Q$. Let $A_t:=A_t(\cE\otimes \bW)$ be the rescaled super connection of the geometric family $\cE\otimes \bW$. With the additional term we must consider the super connection  $A_{\bar Q,t}:=A_t+t^{\frac{1}{2}}\bar Q$. Note that (see \cite[Prop. 10.15]{MR1215720})
\begin{equation}\label{iqwdodqwdwwdqwd1}A_{\bar Q,t}=t^{\frac{1}{2}}(D(\cE\otimes \bW)+\bar Q)+\nabla^{H(\cE\otimes \bW)}-\frac{1}{4t^{\frac{1}{2}}}c(T)\ ,\end{equation}
where $T$ is the curvature tensor associated to the horizontal distribution $T^h\pi$,
and $\nabla^{H(\cE\otimes \bW)}$ is an unitary connection on the Hilbert bundle $H(\cE\otimes \bW)$.

We now calculate the square
\begin{equation}\label{iqwdodqwdwwdqwd}
A_{\bar Q,t}^2=A_t^2+t([D(\cE\otimes \bW),\bar Q]+\bar Q^2)+t^{\frac{1}{2}}[\nabla^{H(\cE\otimes \bW)},\bar Q]\ .
\end{equation}
Note that for two odd endomorphisms $X,Y$ we have by definition $[X,Y]:=XY+YX$.
The mixed term with $\bar Q$ and the curvature $T$ disappears because of $[c(T),\bar Q]=0$ since the Clifford multiplication anti-commutes with $\bar Q$.
We now calculate the commutator terms using a vertical orthonormal frame
$(e_i)$ as well as a horizontal frame $(f_{\alpha})$ and its dual $(f^\alpha)$.
Note that $$\nabla^{H(\cE\otimes \bW)}=\sum_\alpha f^\alpha (\nabla^{S(T^v\pi)\otimes W}_{f_\alpha}+\frac{1}{2}k(f_\alpha))\ ,$$ where $k$ is the mean curvature of the fibre in the direction $f_\alpha$.
Then we have
$$[D(\cE\otimes \bW),\bar Q]=\sum_i[c(e_i)\nabla^{S(T^v\pi)\otimes W}_{e_i},\bar Q]=\sum_i
c(e_i) \nabla^{W}_{e_i}\bar Q$$
and
$$[\nabla^{H(\cE\otimes \bW)},\bar Q]=\sum_\alpha [f^\alpha (\nabla^{S(T^v\pi)\otimes W}_{f_\alpha}
+k(f_\alpha),\bar Q]=\sum_\alpha f^\alpha \nabla^{W}_{f_\alpha}\bar Q\ .$$
The mean curvature drops out since $[k(f_\alpha),\bar Q]=0$ and $[f^\alpha,\bar Q]=0$.
We now perform the Getzler rescaling as in \cite[Ch 10]{MR1215720}.
We only study the terms involving $\bar Q$. We have
\begin{eqnarray}
\lim_{u\to 0}\delta_u (t\bar Q^2)&=&0\label{zzuzuqwdqwdwqd}\\
\lim_{u\to 0}\delta_u(t([D(\cE\otimes \bW),\bar Q])&=&0\nonumber\\
\lim_{u\to 0}\delta_u(t^{\frac{1}{2}}[\nabla^{H(\cE\otimes \bW)},\bar Q])&=&\sum_\alpha f^\alpha  \nabla^{W}_{f_\alpha}\bar Q=:\nabla^h \bar Q\ .
\nonumber\end{eqnarray}
In particular  the limit $\lim_{u\to 0} \delta_u(A_{\bar Q,t}^2)$ exists.
Therefore the $t\to 0$-asymptotic of $\Tr \exp(A^2_{\bar Q,t})$ is still regular.
But we may get a contribution to the local index form
$$\Omega(\cE\otimes \bW,\bar Q):=\varphi\lim_{t\to 0}\Tr  \exp(-A^2_{\bar Q,t})\ ,$$
where $\varphi$ scales $2k$-forms by $\frac{1}{(2\pi i)^k}$.
In formula \cite[10.28]{MR1215720} one has to replace the twisting curvature
$F=R^{\nabla^{W}}$ by $R^{\nabla^W}+ \nabla^h \bar Q$. Note that
$\nabla^W\bar Q\in \Omega^1(\End(W)^{odd})$ commutes 
with multiplication by two-forms, but not necessarily with  $R^{\nabla^{W}}$.
We get
$$ \Omega(\cE\otimes \bW,\bar Q)=\varphi \int_{E/B}  \det^{\frac{1}{2}}\left( \frac{\frac{R^{\nabla^{T^v\pi}}}{2}}{\sinh(\frac{R^{\nabla^{T^v\pi}}}{2})}\right) \tr \exp\left(-R^{\nabla^W}- \nabla^h\bar Q\right)\ .$$
We calculate the $4$-form component of the integrand. Note that $\dim W=\dim W^+-\dim W^-=0$.
\begin{eqnarray*}\lefteqn{
\left[ \det^{\frac{1}{2}}\left( \frac{\frac{R^{\nabla^{T^v\pi}}}{2}}{\sinh(\frac{R^{\nabla^{T^v\pi}}}{2})}\right)  \tr \exp\left(-R^{\nabla^W}\right)\right]_{4}}&&\\&=&  \left[\tr \exp\left(-R^{\nabla^W}- \nabla^h\bar Q\right)\right]_{4}\\
&=& \frac{1}{2}\tr(R^{\nabla^W})^2- \frac{1}{2}\tr(R^{\nabla^W}(\nabla^h \bar Q)^2)+\frac{1}{24} \tr \nabla^h (\bar Q)^4 \ ,
\end{eqnarray*}
where we use that
$\tr(\nabla^h\bar Q R^{\nabla^W} \nabla^h\bar Q)=\tr(R^{\nabla^W} (\nabla^h\bar Q)^2)$.
Note that $$\Omega^2(\cE\otimes \bW)=\frac{1}{2\pi i}\int_{E/B}  \frac{1}{2}\tr(R^{\nabla^W})^2$$ 
and the form $\tr \nabla^h (\bar Q)^4$ has no vertical component
so that
\begin{eqnarray*}
\Omega^2(\cE\otimes \bW,\bar Q)=  \Omega^2(\cE\otimes \bW)-\frac{1}{4\pi i} \int_{E/B} \tr(R^{\nabla^W}(\nabla^h \bar Q)^2)   \ .
\end{eqnarray*}
This gives the curvature of the Pfaffian and determinant line bundle of the family of generalised Dirac operators
\begin{eqnarray}
\frac{1}{2\pi i}R(\Pfaff(\cE\otimes \bW,\bar Q,J))&=&-\frac{1}{2}\Omega^2(\cE\otimes \bW)+\frac{1}{8 \pi i}\int_{E/B} \tr(R^{\nabla^W}(\nabla^h \bar Q)^2) \label{udiwdqwdwqd333}\\\frac{1}{2\pi i}R(\det(\cE\otimes \bW,\bar Q))&=&\Omega^2(\cE\otimes \bW)-\frac{1}{4 \pi i}\int_{E/B} \tr(R^{\nabla^W}(\nabla^h \bar Q)^2) \ .\nonumber\end{eqnarray}

\subsection{Construction of local sections of the Pfaffian}

If $U\subseteq B$ is open  such that restriction of $V$ to $E_U$ is trivial as a spin bundle as in Section \ref{ll} we let $I_U$ denote the set of trivialisations. For every
$Q\in I_U$ we are going to construct a section $d_Q\in C^\infty(U,\Pfaff(\cE\otimes \bW,J))$.
Let $\pr:\R\times U\to U$ be the projection and consider the family
$\pr^*(\cE\otimes \bW)$ over $\R\times U$. Its underlying bundle is $\R\times E_U\to \R\times U$, and we consider the projection $\Pr:\R\times E_U\to E_U$ of total spaces.
  We
 define the odd selfadjoint endomorphism
$\tilde  Q\in \End(\Pr^* W_{E_U})$ so that it equals $a \bar Q$ on the slice
$\{a\}\times E_U\subset \R\times E_U$.
We calculate the Laplacian $\Delta((a,u))$ by extracting the zero-form part of  (\ref{iqwdodqwdwwdqwd}) over the base point $(a,u)\in \R\times U$ at time $t=1$.
We get
$$\Delta((a,u))=D(\cE\otimes \bW)^2+a c(\nabla^W\bar Q)+a^2 \bar Q^2\ .$$
Note that $ \bar Q^2=1$ is positive. For large $a$ the term $a^2 \bar Q^2$ dominates $a c(\nabla^W\bar Q)$. More precisely,
if we assume that $U$ has a compact closure in $B$, then there exists $a_0\ge 0$ such that for 
$a_0\le  a$ the operator $\Delta((a,u))$ is positive, and hence invertible.
Therefore we have the section
$$s^{1/2}_{can}\in C^\infty([a_0,\infty)\times U,\Pfaff(\pr^*(\cE\otimes \bW),\tilde Q,\Pr^*J))\ .$$
The norm of $s^{1/2}_{can}$ is given by (\ref{uqwidqwdwqdwd666}), and we consider the unit-norm section
$s^{1/2,0}_{can}:=\|s^{1/2}_{can}\|^{-1}s^{1/2}_{can}$.
We have a canonical identification
$$\Pfaff(\pr^*(\cE\otimes \bW),\tilde Q,\Pr^*J)_{\{0\}\times U}\cong \Pfaff(\cE\otimes \bW,J)_U\ .$$
For $b\ge a_0$ we define the unit-norm section
$$d(b)\in C^\infty(\R\times U,\Pfaff(\pr^*(\cE\otimes \bW),\tilde Q,\Pr^*J))$$ such that $d(b)(a,u)$ is  the parallel transport of
$s_{can}^{1/2,0}(b,u)$ along the path $[0,1]\mapsto ((1-t)b+ta,u)$.
We define $d_Q(b)\in C^\infty(U, \Pfaff(\cE\otimes \bW,J))$ by evaluation of $d(b)$ at $u=0$, i.e.
$d_Q(b)(u):=d(b)(0,u)$ 

We now consider the $\eta^1$-form
$$\eta^1:=\eta^1(\pr^*(\cE\otimes \bW)_{t_{\tilde Q}})\in \Omega^1([a_0,\infty)\times U)\ .$$ Recall  from (\ref{udduqwdqwdqwd777}) that
\begin{equation}\label{uidqdqdwqdwd}
-\pi i \eta^1=\nabla^{\Pfaff(\pr^*(\cE\otimes \bW),\tilde Q,\Pr^*J)}\log s_{can}^{1/2,0}\ .
\end{equation}
We write
$$\eta^1=a^2(da \theta+\lambda)\ ,$$ where
$\theta\in C^\infty([a_0,\infty)\times U)$ and $\lambda\in C^\infty([a_0,\infty)\times U,\pr^*T^*B)$. 
\begin{prop}\label{p1}
For $a\to \infty$ there are asymptotic expansions 
$$\theta=\sum_{n\ge 0} a^{-n} \theta_{-n}\ , \quad \lambda=\sum_{n\ge 0} a^{-n} \lambda_{-n}\ ,$$
where $\theta_i\in C^\infty(U)$, $\lambda_i\in \Omega^1(U)$.
Moreover,
$\theta_0$ and $\theta_{-3}$ are constant.
\end{prop}
\proof
The second assertion will be shown as a consequence of the first
in the proof of Proposition \ref{p2}.
The existence of the asymptotic expansion will be shown later in Subsection \ref{zzuuzuzee}.
We define
$$\tilde d_Q(b):=d_Q(b)\exp\left(i\pi(\frac{1}{3}b^3\theta_{0}+\frac{1}{2}b^2\theta_{-1}+b\theta_{-2}+\log (b) \theta_{-3})\right)\ .$$

\begin{prop}\label{p2}
The limit 
$$d_Q:=\lim_{b\to \infty} \tilde d_Q(b) $$
exists in the $C^1_{loc}$-sense. The connection one-form of the limit is given by 
$$ \nabla^{ \Pfaff(\cE\otimes \bW,J)}\log d_Q=-\pi i\int_{E_U/U}\eta^3(\cW_{t_Q})\ .$$
\end{prop}
\proof
Note that by (\ref{uidqdqdwqdwd}) for $b^\prime\ge b$ we have
$$\frac{d_Q(b^\prime)}{d_Q(b)}=\frac{d(b^\prime)_{|\{b\}\times U}}{d(b)_{|\{b\}\times U}}=\exp\left(-i\pi \int_{[b,b^\prime]\times U/U} \eta^1\right) \ .$$
If we insert the asymptotic expansion of $\eta^1$
we get
$$ \int_{[b,b^\prime]\times U/U} \eta^1=\frac{b^{\prime 3}-b^3}{3} \theta_0+\frac{b^{\prime2}-b^2}{2} \theta_{-1}+(b^\prime-b)\theta_{-2}+(\log(b^\prime)-\log(b))\theta_{-3}+O(b^{\prime -1},b^{-1})\ .$$
It follows that
$$\frac{\tilde d_Q(b^\prime)}{\tilde d_Q(b)}=\exp(i O(b^{\prime -1},b^{-1}))\ .$$
This implies the existence of $\lim_{b\to \infty} \tilde d_Q(b)$.

Now we consider the connection one-forms.
Let $R:=\frac{1}{2\pi i}R^{\Pfaff(\pr^*(\cE\otimes \bW),\tilde Q,\Pr^*J)}$.
Note that  by (\ref{udiwdqwdwqd333}) we have
$$ R=-\frac{1}{2}\pr^* \Omega(\cE\otimes \bW)+ \frac{1}{8 \pi i} \int_{[0,1]\times E_U/[0,1]\times U} \tr (\pr^* R^{\nabla^W} (\nabla^h \tilde Q)^2)\ .$$
Since $\tilde Q$ contains the variable $a$ linearly we see that
$$-2R=ada \wedge \pr^*S+a^2\pr^*T+ \pr^* \Omega(\cE\otimes \bW)$$ for some 
 $S\in \Omega^1(U)$ and $T\in \Omega^2(U)$. 
On $[a_0,\infty)\times U$ we have the identity
$$d\eta^1=-2R\ .$$
If we assume that $\eta^1$ has an asymptotic expansion
as stated in Proposition \ref{p1}, then we get
$$-\sum_{n\ge 0} a^{2-n} da \wedge d\theta_{-n} +\sum_{n\ge 0} a^{2-n}d\lambda_{-n} +\sum_{n\ge 0} (2-n)a^{1-n} da\wedge \lambda_{-n}=ada \wedge \pr^*S+a^2\pr^*T+\pr^* \Omega(\cE\otimes \bW)\ .$$
This gives the following identities for the terms in the asymptotic expansion of $\eta^1$:
 \begin{enumerate}
\item[$a^2$:] $d\lambda_0=\pr^*T$, \hspace{1cm}$d\theta_0=0$
\item[$a^1$:]  $d\lambda_{-1}=0$, \hspace{1cm}$-d\theta_{-1}+2\lambda_0=\pr^*S$
\item[$a^0$:] $d\lambda_{-2}= \pr^* \Omega(\cE\otimes \bW)$,\hspace{1cm} $-d\theta_{-2}+\lambda_{-1}=0$
\item[$a^{-1}$:] $d\lambda_{-3}=0$, \hspace{1cm}$-d\theta_{-3}=0$\ .
\end{enumerate}
In particular we obtain the second assertion of Proposition \ref{p1}.
Let $$\omega(b):=\nabla^{\Pfaff(\pr^*(\cE\otimes \bW),\tilde Q,\Pr^*J)}\log d(b)$$ be the connection one-form of the section $d(b)$.
Since the section $d(b)$ is parallel in the $a$-direction we get for a vector field $X\in \cX(U)$  that
\begin{eqnarray*}
\partial_a \omega(b)(X)d(b)&=&\nabla^{\Pfaff(\pr^*(\cE\otimes \bW),\tilde Q,\Pr^*J)}_{\partial_a} \nabla^{\Pfaff(\pr^*(\cE\otimes \bW),\tilde Q,\Pr^*J)}_X d(b)\\&=&R^{\Pfaff(\pr^*(\cE\otimes \bW),\tilde Q,\Pr^*J)}(\partial_a,X)d(b)\\&=& - \pi ia S(X) d(b)\end{eqnarray*}  and therefore
$$\partial_a \omega(b)(X)=  -\pi ia S(X)\ .$$
Let $\omega_Q(b):=\nabla^{\Pfaff(\cE\otimes \bW,J)} \log d_Q(b)$.
Then we get by integration from $0$ to $b$  
$$ \omega(b)_{\{b\}\times U}(X)=-\frac{\pi ib^2}{2}S(X)+\omega_Q(b)(X)\ .$$
Note that
$$ \omega(b)_{\{b\}\times U}(X)=-\pi i \eta^1_{\{b\}\times U}(X)=-\pi ib^2\lambda(b)(X)\ .$$
The connection one-form $\tilde \omega_Q(b)$ of $\tilde d_Q(b)$ is given by 
\begin{eqnarray*}
\tilde \omega_Q(b)&=&\omega_Q(b)+i\pi\left(\frac{1}{3} b^3 d\theta_0+\frac{1}{2} b^2 d\theta_{-1}+b d\theta_{-2}+ \log b d\theta_{-3}\right)
\\&=&\omega_Q(b)+i\pi\left(\frac{1}{2} b^2 d\theta_{-1}+b d\theta_{-2}\right)\\
&=&-\pi i b^2 \lambda(b)+\frac{\pi i}{2}b^2 S+i\pi\left(\frac{1}{2} b^2 d\theta_{-1}+b d\theta_{-2}\right)\\
&=&\pi i\left(b^2(\frac{1}{2}S-\lambda_0+\frac{1}{2}d\theta_{-1})+b(-\lambda_{-1}+d\theta_{-2})-\lambda_{-2}+o(b^{-1})\right)\\
&=&-\pi i \lambda_{-2} +o(b^{-1})\ . 
\end{eqnarray*}
Therefore
$$\lim_{b\to \infty}\tilde \omega_Q(b)=-\pi i  \lambda_{-2}\ .$$
It remains to calculate
$ \lambda_{-2}$.
To this end we consider a function
$\chi\in C^\infty(0,\infty)$ such that $\chi(t)=0$ for $t\le 1$ and
$\chi(t)=1$ for $t\ge 2$. Let $(x,a,u)\in [0,1]\times \R\times U$ and $\tilde \pr:
[0,1]\times \R\times U\to U$ be the projection. Then
$\tilde \pr^*E_U\cong  [0,1]\times \R\times E_U$. We let $\tilde \Pr: [0,1]\times \R\times E_U\to E_U$ be the projection. On  $H(\tilde \pr^*(\cE\otimes \bW))$  we consider the family of rescaled super connections
$\tilde A_t$ which is given by
\begin{equation}\label{ueiqweuqwe}
\tilde A_t:=\tilde \pr^* A_t(\cE\otimes \bW)+ t^{\frac{1}{2}}
a(x+ (1-x)\chi( t^{\frac{1}{2}}a))\tilde \Pr^* \bar Q\ .
\end{equation}
The local index theory and the proof of  Proposition \ref{p1} given in Subsection \ref{zzuuzuzee} still applies to this more general super connection. The additional term involving the cut-off function $\chi$ vanishes
in the Getzler rescaling (\ref{zzuzuqwdqwdwqd}) and does not contribute to the curvature.
We want to show that $\lambda_{-2}$ is independent of $x$.
Note that as $a\to \infty$ we have an asymptotic expansion
$$\tilde \eta^1:=\eta^1(\tilde A_t)\sim a^2\sum_{n\ge 0} a^n(da\wedge \tilde \theta_{-n} 
+ dx\wedge \tilde \kappa_{-n}+ \tilde \lambda_{-n})\ ,$$
where the terms may depend on $x$.
We have
$$d\tilde \eta^1=(ax^2 da+a^2x dx) \wedge \tilde \pr^*S+a^2x^2 \tilde \pr^*T+\tilde \pr^*\Omega(\cE\otimes \bW)\ .$$
From this we deduce (note that the term $\tilde \kappa_{-2}$ does not contribute) that
$$\partial_x \tilde \lambda_{-2}=0\ .$$

The restriction of $\tilde A_t$ to the slice $\{x=0\}$ is the super connection of a tamed
geometric family, and the parameter enters  as appropriate for adiabatic limits, see \cite[2.2.5]{bunke-20071}.
We get
$$\lambda_{-2}=\tilde \lambda_{-2|\{1\}\times U}=\tilde \lambda_{-2|\{0\}\times U}=\lim_{a\to \infty} \tilde \eta^1_{|\{(0,a)\}\times U}=\int_{E_U/U} \eta^3(\cW_{t_Q})\ .$$
\hB

We now consider a second trivialisation $Q^\prime\in I(U)$ and define the section  $d_{Q^\prime}$ by Proposition \ref{p2}. We assume that $U$ is contractible. Then we can choose a homotopy
$H$ from $Q$ to $Q^\prime$. It induces a taming of the family
$\pr^*(\cE\otimes \bW)$, where $\pr:[0,1]\times U\to U$ is the projection. It furthermore induces a taming $\Pr^*\cW_{t_H}$, where
$\Pr:[0,1]\times E_{U}\to E_{U}$ is the induced projection.

\begin{lem}
We have
$$\frac{d_{Q^\prime}}{d_Q}=\exp\left(-\pi i \int_{[0,1]\times E_{U}/U} \eta^3(\Pr^*\cW_{t_H})\right)\ .$$
\end{lem}
\proof
Indeed we can define the section $d_H\in C^\infty([0,1]\times U, \Pfaff(\pr^*(\cE\otimes \bW),\Pr^*J))$ by Proposition \ref{p2}.
 It restricts to $d_Q$ and $d_{Q^\prime}$ at
$\{0\}\times U$ and $\{1\}\times U$.
Of course,
$$\Pfaff(\pr^*(\cE\otimes \bW),\Pr^*J)\cong \pr^* \Pfaff(\cE\otimes \bW,J)\ .$$
Therefore
$$\frac{d_{Q^\prime}}{d_Q}=\exp\left(\int_{[0,1]\times U/U} \nabla^{\Pfaff(\pr^*(\cE\otimes \bW),\Pr^*J)} \log d_H\right)$$
But again by Proposition \ref{p2} we have
$$\nabla^{\Pfaff(\pr^*(\cE\otimes \bW),\Pr^*J)} \log d_H= -\pi i \int_{[0,1]\times E_{U}/[0,1]\times U}  \eta^3(\Pr^*\cW_{t_H})\ .$$
This gives the result.
\hB 

\subsection{Proof of Theorem \ref{hjjhdasd}}

Recall the construction of the geometric line bundle $\bL$ in Section \ref{ll}.
\begin{theorem}\label{mmm1}
There exists a canonical functorial isomorphism of geometric line bundles
$$\bL\cong \Pfaff(\cE\otimes \bW,J)\ .$$
It is characterised by the property that for every open $U\subseteq B$ and $Q\in I_U$ it maps the section 
$s_Q\in C^\infty(U,L)$ to the section $d_Q\in C^\infty(U,\Pfaff(\cE\otimes \bW,J))$.
\end{theorem}
\proof
By inspection one checks that the cocycles 
and the connection one-forms for the  collections of local sections
$(s_Q)_{U,Q\in I_U}$ and $(d_Q)_{U,Q\in I_U}$ coincide.
Note that all constructions are natural with respect to pull-back along smooth maps $B^\prime\to B$.
\hB

%

\subsection{Asymptotic expansion of $\eta$-forms in the adiabatic limit}\label{zzuuzuzee}
In this technical subsection we prove the asymptotic expansions of the $\eta^1$-form stated in Proposition \ref{p1} and used in the more general case in the course of the proof of Proposition \ref{p2}.
In general the $\eta^1$-form for a rescaled super connection $\tilde A_t$ is given by 
\begin{equation}\label{iwioqdqwdqwdqwd773}
 \eta^1:=\eta^1(\tilde A_t)   =\frac{1}{2\pi i}\int_0^\infty\Tr \left[\partial_t \tilde A_t \ee^{-\tilde A_t^2 }\right]_1\ .
\end{equation}
In our situation we consider (\ref{ueiqweuqwe}), i.e. 
$$\tilde A_t:=\tilde \pr^* A_t(\cE\otimes \bW)+ t^{\frac{1}{2}}
a(x+ (1-x)\chi( t^{\frac{1}{2}}a))\tilde \Pr^* \bar Q\ .$$
We simplify the notation. We write
$$\tilde A_t:=A_t+ f(t^{\frac{1}{2}}a)  \bar Q\ ,$$
where $A_t$ is a Bismut super connection associated to a family of Dirac operators $D$ over a base $\R\times B$ with coordinates $(a,b)$. We let $H:=\nabla^h\bar Q$ be the horizontal derivative in the $B$-direction of $\bar Q$, where  $\bar Q$ is an odd involution, and $f\in C^\infty([0,\infty)\times B)$
is such that
$f(t,b)=t$ for $t\ge 1$.
Then we have for the one-form component
\begin{eqnarray*}
[\Tr\partial_t \tilde A_t \ee^{-\tilde A_t^2} ]_1&=&\left[\Tr
\frac{1}{2t^{\frac{1}{2}}} (D+af^\prime(t^{\frac{1}{2}}a)\bar Q)\ee^{-
(t^{\frac{1}{2}}D+f(t^{\frac{1}{2}}a)  \bar Q)^2 -t^{\frac{1}{2}}f (t^{\frac{1}{2}}a)H-t^{\frac{1}{2}}f^\prime(t^{\frac{1}{2}}a)\bar Q da
}\right]_1\\
&=&-\frac{1}{2}\Tr
(D+af^\prime(t^{\frac{1}{2}}a)\bar Q)\ee^{-
(t^{\frac{1}{2}}D+f(t^{\frac{1}{2}}a)  \bar Q)^2} \left\{  f (t^{\frac{1}{2}}a)H+ f^\prime(t^{\frac{1}{2}}a)\bar Q da\right\}
\end{eqnarray*}
We further calculate

$$(t^{\frac{1}{2}}D+f(t^{\frac{1}{2}}a)  \bar Q)^2=tD^2+t^{\frac{1}{2}}f(t^{\frac{1}{2}}a)E+f(t^{\frac{1}{2}}a)^2\ ,$$
where $E:=c(\nabla^v\bar Q)$ is the Clifford multiplication by the vertical
derivative of $\bar Q$. 

We now introduce the variable $s=a^{\frac{3}{2}}t$.
Then we have $dt=a^{-\frac{3}{2}}ds$, $t^{\frac{1}{2}}a=a^{\frac{1}{4}}s^{\frac{1}{2}}$, and 
\begin{eqnarray*}
\eta^1&=&-\frac{1}{4\pi i}\int_0^\infty 
\Tr\left\{(D+af^{\prime}(a^{\frac{1}{4}}s^{\frac{1}{2}}) \bar Q)
\ee^{-a^{-\frac{3}{2}}sD^2+a^{-\frac{3}{4}}s^{\frac{1}{2}} f(a^{\frac{1}{4}}s^{\frac{1}{2}})E+f(a^{\frac{1}{4}}s^{\frac{1}{2}})^2}\right.\\&&\left.(
f(a^{\frac{1}{4}}s^{\frac{1}{2}}) H+f^\prime(a^{\frac{1}{4}}s^{\frac{1}{2}})\bar Q da)\right\} a^{-\frac{3}{2}}ds\ .
\end{eqnarray*}
Note that for $s\ge 1$ and $a\ge 1$ we have
$f(a^{\frac{1}{4}}s^{\frac{1}{2}})=a^{\frac{1}{4}}s^{\frac{1}{2}}$ and $f^\prime(a^{\frac{1}{4}}s^{\frac{1}{2}})=1$.
In this region the integrand simplifies to 
$$\ee^{-a^{\frac{1}{2}}s}\Tr\left\{(D+a \bar Q)
\ee^{-sa^{-\frac{3}{2}}(D^2+aE)} 
(a^{\frac{1}{4}} s^{\frac{1}{2}} H+\bar Q da)\right\} a^{-\frac{3}{2}}\ .$$

\begin{lem}\label{zwduqwdwqdwqdwqdwd}Locally uniformly in $B$
there exists constants $C<\infty$ and $c>0$ independent of $a$ such that
$$|\int_{1}^\infty \dots ds|\le C e^{-c a^{\frac{1}{2}}}\ .$$
\end{lem}
\proof
We will show that
there exists $C<\infty$ and $M\in \nat$ such that for all $s\ge 1$ and $a\ge 1$
$$|\Tr\left\{(D+a \bar Q)
\ee^{-sa^{-\frac{3}{2}}(D^2+aE)} 
(a^{\frac{1}{4}} s^{\frac{1}{2}} H+\bar Q da)\right\}|\le C a^Ms^M\ .$$
The growth of the right-hand side in this estimate can be absorbed in the prefactor
$\ee^{-a^{\frac{1}{2}}s}$. The assertion of the Lemma then follows by elementary calculus.

This estimate of the trace is obtained by a combination of spectral and trace class estimates. The inclusion of  Sobolev spaces $H^k\to H^0$ is of trace class if $k$ is larger then the  dimension of the underlying space which is two in our case. 
We further use ellipticity of $D $ which implies that the graph norm associated to $D^k$ is equivalent to the $k$'th Sobolev norm.
We write $A:=D+aE$. It suffices to show that for every $N\in \nat $ there exist constants
$c>0$ and $C<\infty$ and integers $M\in \nat$
\begin{equation}\label{ztqszuqqwdqwdqw}
\|D^{2N}\ee^{-s a^{-\frac{3}{2}}A}\|\le C  s^Ma^M \ .
\end{equation}

We use $u:=s a^{-\frac{3}{2}}$ and write
\begin{equation}\label{zuiudwqdqwdwqd}
D^{2N}e^{-uA^2}=A^{2N}e^{-uA^2}+(D^{2N}-A^{2N})e^{-uA^2}\ .
\end{equation}
The first summand can be written as
$$u^{-N} (u^{\frac{1}{2}}A)^{2N}e^{-(u^{\frac{1}{2}}A)^2}\ .$$
We now use that
$$\|(u^{\frac{1}{2}}A)^{2N}e^{-(u^{\frac{1}{2}}A)^2}\|\le C $$
This follows by the spectral mapping principle from the bound $$\sup_{x\in [0,\infty]}x^{2N}e^{-x^2}<\infty\ .$$  
If we combine these estimates we get
$$\|A^{2N}e^{-uA^2}\|\le  u^N C\ .$$
Note that the difference
$(D^{2N}-A^{2N})$ can be written as
$\sum_{i=0}^{2N-2} a^{k_i}E_i$, 
where $E_i$ are differential operators of order $i$. 
We can therefore use induction by $N$ in order to deal with the second term in (\ref{zuiudwqdqwdwqd}) and to get an estimate of the form
$$\|D^{2N} e^{-uA^2}\|\le C u^M a^{M} $$
for some $M\in \nat$. We now insert $u=sa^{-\frac{3}{2}}$ and get 
(\ref{ztqszuqqwdqwdqw}).
\hB

\begin{lem} Locally uniformly in $B$ 
we have an asymptotic expansion
$$|\int_0^{1}\dots ds|\sim  a^2\sum_{n\ge 0} a^{-n} \eta^1_{-n}\ .$$
\end{lem}
\proof
We must control the integral kernel of the smoothing operator
$$\ee^{-t(D^2+t^{-\frac{1}{2}} f(t^{\frac{1}{2}}a) E)}$$
in the region $0< t\le a^{-\frac{3}{2}}$.  
For fixed $T$ we  construct  a formal solution of the heat equation
$$(\partial_t-(D^2+T  E))H_t=0\ .$$
by the iterative procedure of the proof of \cite[Thm 2.26]{MR1215720} and keep control of the dependence on $T$.
We get
$$H(t,x,y)=q_t(x,y)\sum_{n\ge 0} t^i \Phi_i(x,y)\ ,$$
where $q_t$ is a Gauss kernel and
the coefficients $ \Phi_i(x,y)$ implicitly also depend on $T$.
The coefficients  $\Phi_i(x,y)$ are given by an iterative formula
stated in  \cite[Thm 2.26]{MR1215720}. By inspection we see that
it is a polynomial in $T$ of degree at most $i$.
If we write
$$(\partial_t-(D^2+T E))H_t=q_t(x,y)\sum_{n=0}^N t^i \Phi_i(x,y)+ t^{N-1}r_t^N(x,y)\ ,$$
then the remainder term is bounded in $C^l$-norm by $t^{-l}$.
Moreover, it is a polynomial in $T$ of degree at most $N+1$.
Note that the $t$-power is explained by $N-1=N-\frac{\dim(E/B)}{2}$.
We now construct the heat kernel
$\ee^{-t(D^2+T  E)}$ using the Volterra series method as in 
 \cite[Sec. 2.4]{MR1215720}. Then we set
$T(t,a):=t^{-\frac{1}{2}} f(t^{\frac{1}{2}}a)$.
We observe that for $N>\frac{M}{2}$
\begin{equation}\label{uiqwdqwdqwdwqdqwd}
\sup_{t\in (0,a^{-\frac{3}{2}})}t^N T^M(t,a)\le C a^{-\frac{3N}{2}+M}\ .
\end{equation}
We split
$$t^{N-1}r^N_t=\left(t^{\frac{N-1}{4}} \right)\left( t^{\frac{3N-3}{4}} r^N_t\right)\ .$$
Using that $r^N_t$ is a polynomial in $T$ of degree $N+1$ we apply (\ref{uiqwdqwdqwdwqdqwd}) to the second factor  and see that for $N>3$ the remainder term estimate in
\cite[Thm. 2.20 (3)]{MR1215720} is for $t\in (0,a^{-\frac{3}{2}})$
$$\|r_t^N\|_l\le C(l)  a^{3-\frac{N}{8}} t^{\frac{N-1}{4}-\frac{l}{2}}\ ,$$
where $C(l)$ does not depend on $a$ furthermore.
This leads to an estimate
$$\|r^{k+1}_t\|_l \le C^{k+1} a^{k(3-\frac{N}{8})}
t^{(k+1)(\frac{N-1}{4}-\frac{l}{2})} \frac{t^k}{k!}$$
in \cite[Lemma 2.21]{MR1215720}.
Finally, in \cite[Theorem 2.23 (2)]{MR1215720} we get for all integers $l,n$, that
 for sufficiently large $N$ depending on $n,l$ the  approximate kernel $$k^N_t(x,y):=q_t(x,y)\sum_{n=0}^N t^i \Phi_i(x,y)$$
differs from the true kernel in $C^l$-norm by 
$Ca^{-n}$ uniformly in $t$ and $a$.

In order to get the asymptotic expansion of 
the $\eta^1$-form we can therefore replace the true heat kernel by its 
approximation.  We must derive an asymptotic expansion of
\begin{eqnarray*}
\eta^1&\sim&-\frac{1}{4\pi i}\int_0^1 
\Tr\left\{(D+af^{\prime}(a^{\frac{1}{4}}s^{\frac{1}{2}}) \bar Q)
k^N_{a^{-\frac{3}{2}}s}\ee^{ -f(a^{\frac{1}{4}}s^{\frac{1}{2}})^2}\right.\\&&\left.(
f(a^{\frac{1}{4}}s^{\frac{1}{2}}) H+f^\prime(a^{\frac{1}{4}}s^{\frac{1}{2}})\bar Q da)\right\} a^{-\frac{3}{2}}ds\\&=&
-\frac{1}{4\pi i}\int_0^{a^{\frac{1}{4}}} 
\Tr\left\{(D+af^{\prime}(u) \bar Q)
k^N_{a^{-2}u^2}\ee^{ -f(u)^2}\right.\\&&\left.(
f(u) H+f^\prime(u)\bar Q da)\right\} a^{-2} 2udu\\&=&
-\frac{1}{4\pi i}\int_0^{a} 
\Tr\left\{(D+af^{\prime}(u) \bar Q)
k^N_{a^{-2}u^2}\ee^{ -f(u)^2}\right.\\&&\left.(
f(u) H+f^\prime(u)\bar Q da)\right\} a^{-2} 2udu\  ,
\end{eqnarray*}
where in the last step we use that because of the 
factor $\ee^{-f(u)}$ the integral $\int_{a^{\frac{1}{4}}}^a\dots du$ does not contribute to the asymptotic expansion.
Also note that $T(t,a)=au^{-1}f(u)$.
Therefore the integrand obviously  has an expansion in terms of powers
of $a$ with integrable functions of $u$\footnote{We a priori know that
all terms are integrable at $u=0$.}. So the integral gives an expansion in powers of $a$.
It remains to determine the leading order. It is a multiple of $a^2$. 
\hB

\section{String structures}\label{str}

\subsection{Geometric and multiplicative gerbes}

In this subsection we recall some aspects of the theory of $U(1)$-banded gerbes in manifolds with an emphasis on geometric structures.
Furthermore, we review multiplicative gerbes on Lie groups in some detail.

We use the language of stacks and  refer to \cite{MR2206877}
for a very readable introduction to stacks in manifolds which suffices for the purpose of the present paper. This in particular applies to the notion of a gerbe
with band in an abelian Lie group.  

In the following a gerbe is a $U(1)$-banded gerbe in manifolds.
Gerbes over a given manifold $M$ form a monoidal $2$-category. Details of the construction of the tensor product can be found in \cite[6.1.9]{bunke-2008-2006}.
The tensor unit is the trivial gerbe $M\times \cB U(1)$ with the quotient stack $\cB U(1):=[*/U(1)]$.

For every pair  $t_1,t_2:\cH\to \cH^\prime$ of $1$-morphisms between gerbes there is an associated  $U(1)$-principal bundle which we denote by $\frac{t_2}{t_1}$. Since this bundle plays an important role in the description of geometric structures let us describe this bundle in greater detail by  characterizing
  its sheaf of sections. By $\cU(1)$ we denote the sheaf of $U(1)$-valued smooth functions.
 For $U\subseteq M$ and an object
$o\in \cH(U)$ in the groupoid $\cH(U)$ we consider the set-valued presheaf $$U\supseteq V\mapsto \underline{\Hom} (t_1(o),t_2(o))(V):=\Hom_{\cH^\prime(V)}(t_1(o)_{|V},t_2(o)_{|V})$$
on $U$.
Since $\cH^\prime$ is a stack
it is a sheaf, and since $\cH^\prime$ is a $U(1)$-banded gerbe, it is in addition 
  is a sheaf of $\cU(1)$-torsors. This sheaf is independent of the choice of $o$ in the sense, that for any other  choice $o^\prime\in \cH(U)$ and isomorphism $\phi\in \Hom_{\cH(U)}(o,o^\prime)$   we have
an isomorphism of sheaves 
$\underline{\Hom} (t_1(o),t_2(o))\cong  \underline{\Hom} (t_1(o^\prime),t_2(o^\prime))$
which is independent of $\phi$. By the axioms for a gerbe, we can cover
$M$ by open subsets $U\subseteq M$ for wich $\cH(U)$ has objects. Moreover, any two objects of $\cH(U)$ are locally isomorphic. 
We can therefore glue the sheaves $\underline{\Hom} (t_1(o),t_2(o))$ to a global sheaf on $M$ which is the sheaf of sections of the $U(1)$-principal bundle $\frac{t_2}{t_1}$.

A $2$-morphism $t_1\Rightarrow t_2$ can be viewed as a trivialisation   $1_M\to \frac{t_2}{t_1}$,
where $1_M:=M\times U(1)$ is the trivial $U(1)$-bundle.

The isomorphism class of a gerbe $\cH$  in  the $2$-category of gerbes over $M$ is classified by the Dixmier-Douady class $DD(\cH)\in H^3(M;\Z)$. We have
$DD(\cH\otimes \cH^\prime)=DD(\cH)+DD(\cH^\prime)$.
Furthermore, if  $\cH$ and $\cH^\prime$ are isomorphic, then the set of isomorphism classes in the category $\Hom(\cH,\cH^\prime)$ is a torsor over $H^2(M;\Z)$.

Next we discuss multiplicative gerbes on a Lie group $G$. The case of interest in the present paper is $G=Spin(n)$.
The notion of a multiplicative gerbe in a simplicial context has been introduced in 
 \cite{MR2174418}. Here we prefer to work directly in the $2$-catgory of $U(1)$-banded gerbes on $G$. 

For $k\ge l$ let  $\pr_{i_1\dots i_l}:G^k\to G^l$ denote
the projection
$(g_1,\dots,g_k)\mapsto (g_{i_1},\dots,g_{i_l})$, and for $1\le i\le k-1$ let $m_{i,i+1}:G^k\to G^{k-1}$  denote the multiplication
$(g_1,\dots,g_k)\mapsto (g_1,\dots,g_ig_{i+1},\dots,g_k)$.
 We use the notation $\cG_{i}:=\pr_i^*\cG$ and $m:=m_{12}$.
\begin{ddd}
A multiplicative structure on a gerbe $\cG$  on $G$ is given by a $1$-morphism
\begin{equation}\label{iudqwdwqdwqd}
\mu:\cG_1\otimes \cG_2\to m^*\cG
\end{equation} and an associativity $2$-morphism satisfying a higher coherence condition.  
\end{ddd} 

In the simplicial context the following Lemma has been shown in \cite[Sec. 5]{MR2174418}. 
\begin{lem}\label{udiqwdqwdwqddqwdqwdqwd}
If $G$ is compact, connected and simply connected, then a gerbe
$\cG$ on $G$ admits a multiplicative structure which is unique up to isomorphism.
\end{lem}
\proof
We first observe, using the K\"unneth formula and $H^2(G;\Z)=0=H^1(G;\Z)$, that every class
$x\in H^3(G;\Z)$ is primitive in the sense that it satisfies
$m^*x=\pr_1^*x+\pr_2^*x$.
Therefore the Dixmier-Douday classes of the two sides in (\ref{iudqwdwqdwqd})
coincide so that we can choose an isomorphism $\mu$. 
Since $H^2(G^2;\Z)=0$ it is unique up to isomorphism.

Next we must find an associativity $2$-morphism.
The associativity of the group multiplication 
can be written as $m\circ m_{12}=m\circ m_{23}$. We have
two $1$-morphisms 
\begin{equation}\label{zuiqwkdqwdqwdqwdqwdd}
\xymatrix{& (m\circ  m_{12})^*\cG\ar@{=}[dd]\\ \cG_1\otimes \cG_2\otimes \cG_3\ar@/^0.5cm/[ur]^{t_0}\ar@/_0.5cm/[dr]_{t_1}&\\
& (m\circ  m_{23})^*\cG}
\end{equation}
given by 
$t_0:=\mu\circ \mu_{12}$ and $t_1:=\mu\circ\mu_{23}$, where $\mu_{12}$ has a meaning  analogous to $m_{12}$.
Since $H^2(G^3;\Z)=0$ the $U(1)$-bundle $\frac{t_1}{t_0}$ is trivializable.
We choose a trivialisation $\tilde a$.
The associativity morphism $a:t_0\Rightarrow t_1$ is a trivialization of this $U(1)$-bundle which satisfies a higher coherence condition. 
We refrain from writing out this diagram but note that the deviation from $\tilde a$ satisfying the higher coherence condition is a smooth group $4$-cocycle
$c_{\tilde a}\in C_{gr}^4(G, U(1))$, where for an abelian  Lie group $A$ we write $C_{gr}^i(G, A):=C^\infty(G^{i},A)$. The differential of the complex
$C_{gr}^*(G, A)$ is given for $c\in C_{gr}^i(G,A)$ by 
$$\delta c(g_1,\dots,g_{i+1}):=c(g_2,\dots,g_{i+1})+\sum_{j=1}^{i}(-1)^j c(g_1,\dots,g_{j}g_{j+1},\dots,g_{i+1})+(-1)^{i+1}c(g_1,\dots,g_{i})\ .$$
The cohomology $H^*_{gr}(G;U(1))$ of the complex 
$(C_{gr}^*(G, A),\delta)$ is the smooth group cohomology of $G$ with coefficients in $U(1)$.  We claim that $H^i(G;U(1))=0$ for $i\ge 1$.
Since $G$ is simply connected, the sequence of coefficients
$$0\to \Z\to \R\to U(1)\to 0$$ induces  a long exact sequence in smooth group cohomology
$$H_{gr}^i(G;\R)\to H_{gr}^i(G;U(1))\to H^{i+1}_{gr}(G;\Z) \to 
H^{i+1}_{gr}(G;\R)\ .$$ 
Since $G$ is compact its higher smooth group cohomology with coefficients the real vector space $\R$ vanishes by the usual averaging argument. 
Furthermore, since $G$ is connected and $\Z$ is discrete, we have
$C_{gr}^*(G,\Z)\cong C_{gr}^*(\{1\};\Z)$ so that 
$H^i_{gr}(G;\Z)\cong H^i_{gr}(\{1\};\Z)=0$ for $i\ge 1$.   This implies the claim.

From $H_{gr}^4(G;U(1))=0$ we deduce
that there is a function $d:G^3\to U(1)$ such that $\delta d=c_{\tilde a}$. If we set $a:=\tilde a d^{-1}$, then $a$ satisfies the required higher coherence relation. This shows the existence of a multiplicative structure.

We now discuss uniqueness. If we choose another multiplication $\mu^\prime$
with associativity isomorphism $a^\prime$, then we get a $U(1)$-bundle $\frac{\mu^\prime}{\mu}$ over $G^2$ which is trivialisable. Let us fix a trivialization $\phi:1_{G^2}\to \frac{\mu^\prime}{\mu}$. Then we can use this $2$-isomorphism $\phi:\mu\Rightarrow \mu^\prime$ in order to compare $a$ with $a^\prime$, say to view $a$ and $a^\prime$ as trivializations of the same bundle.
Then $a^\prime=a  z$ for some cocycle $z\in C^3_{gr}(G,U(1))$.
Since $H^3_{gr}(G;U(1))=0$ we see that $z=\delta v$ for a chain $v\in C^2_{gr}(G;U(1))$. If we modify $\phi$ by $v$, then we get a modified trivialization $\phi^\prime$. If compared with $\phi^\prime$, then $a$ and $a^\prime$ become equal.
\hB 

Below we will apply Lemma \ref{udiqwdqwdwqddqwdqwdqwd} to the group $G=Spin(n)$ for $n\ge 3$ and the basic gerbe $\cG$ whose Dixmier-Douady class is a generator of $H^3(Spin(n);\Z)\cong \Z$.

Geometric structures on $U(1)$-banded gerbes have been popularised in \cite{MR1876068}, \cite{MR1197353} and are usually described in particular models, e.g. of a bundle gerbe.
In the following we give a model-independent account.
By $0_M:=M\times \cB U(1)$ we denote the trivial gerbe on $M$.

\begin{ddd}
A connection $\omega$ on a gerbe $\cH$ over a manifold $M$ associates to every local trivialisation
$t:0_U\stackrel{\sim}{\to} \cH_U$, $U\subseteq M$,  a $2$-form $\omega_t$. This association must be  compatible with restriction to subsets.  Furthermore,  to a pair $(t_0,t_1)$ of local trivialisations of  the gerbe $\cH$ the connection $\omega$ associates  a connection $\nabla^{\frac{t_1}{t_0},\omega}$ on the associated $U(1)$-bundle
$\frac{t_1}{t_0}$ such that $$R^{\nabla^{\frac{t_1}{t_0},\omega}}=\omega_{t_1}-\omega_{t_0}\ .$$ 
This association is again compatible with restriction to open subsets.
\end{ddd}

The form $d\omega_t\in \Omega^3(U)$ is independent of $t$ and therefore the restriction of a closed global three-form
$R^\omega\in \Omega^3(M)$, the curvature of the connection $\omega$. The cohomology class of $R^\omega$  is the image of the Dixmier-Douady class $DD( \cH)$ in de Rham cohomology.

If $\alpha\in \Omega^2(M)$, then we can define a new connection $\omega+\alpha$ such that
$(\omega+\alpha)_t=\omega_t+\alpha$ and $\nabla^{\frac{t_1}{t_0},\omega+\alpha}=\nabla^{\frac{t_1}{t_0},\omega}$. We have $R^{\omega+\alpha}=R^\omega+d\alpha$.

Next we discuss the notion of a connection on a morphism between geometric gerbes.  Such a notion has first been introduced in a slightly different way in \cite{MR2318389}. 
Let $f:\cH\to \cH^\prime$ be a morphism between gerbes with connections
$\omega$ and $\omega^\prime$.

\begin{ddd} A connection $\omega_f$ on the morphism $f$ 
associates a connection $\nabla^{\frac{t^\prime}{f\circ t},\omega_f}$ on the  $U(1)$-bundle
$\frac{t^\prime}{f\circ t}$ for each pair of local trivialisations
$t:0_U\to \cH_U$ and $t^\prime:0_U\to \cH_U^\prime$ such that
 for every other pair $\tilde t,\tilde t^\prime$ of local trivializations
we have
\begin{equation}\label{uifwefwefewff}
(\frac{\tilde t^\prime}{f\circ \tilde t},\nabla^{\frac{\tilde t^\prime}{f\circ \tilde t},\omega_f})\cong 
(\frac{\tilde t^\prime}{t^\prime},\nabla^{\frac{\tilde t^\prime}{t^\prime},\omega^\prime})\otimes  (\frac{t^\prime}{f\circ t},\nabla^{\frac{t^\prime}{f\circ t},\omega_f})\otimes (\frac{t}{\tilde t},\nabla^{\frac{t}{\tilde t},\omega})
\end{equation} as $U(1)$-bundles with connection.
This association must be compatible with restriction to open subsets.
\end{ddd}
By the relation (\ref{uifwefwefewff})
there is a unique  form $R^{\omega_f}\in \Omega^2(M)$ such that
$$R^{\omega_f}_U= \omega^\prime_{t^\prime}-\omega_{t}-R^{\nabla^{\frac{t^\prime}{f\circ t},\omega_f}} $$
for all pairs of objects $t\in \cH(U)$ and $t^\prime\in \cH^\prime(U)$.
This form is called
 the  curvature of the connection.

Note that connections always exist. The curvature satisfies
\begin{equation}\label{udiqwdqwdwqdwqd}
dR^{\omega_f}=R^{\omega^\prime}-R^{\omega}\ .
\end{equation}
Two gerbes $(\cH,\omega)$, $\cH^\prime,\omega^\prime)$ with connection are
isomorphic, if there exists a morphism
$f:\cH\to \cH^\prime$ which admits  connection
$\omega_f$ with the following two properties
\begin{enumerate}
\item $\omega_f$ is flat, and 
\item\label{uiduwqidhqwdwqdwqdqd} for each  $t$ it associates to the pair $(t,f\circ t)$ 
  a trivial bundle with connection $(\frac{f\circ t }{f\circ t},\nabla^{\frac{ f\circ t}{f\circ t},\omega_f})$.
\end{enumerate}
Note that by (\ref{uifwefwefewff})
this  condition fixes the connection $\omega_f$ uniquely.
In  \cite[Def. 4.2.3]{MR2318389} such a connection is called
compatible.

If $f_0,f_1:\cH\to \cH^\prime$ are two $1$-morphisms with connections $\omega_{f_0}$ and $\omega_{f_1}$, then the bundle
$\frac{f_1}{f_0}$ has an induced connection
$\nabla^{\frac{f_1}{f_0},\omega_{f_1},\omega_{f_0}}$ with curvature $R^{\omega_{f_0}}-R^{\omega_{f_1}}$.
This shows that the curvature of a connection on a morphism
between two gerbes with connections is determined up to exact forms by (\ref{udiqwdqwdwqdwqd}).
There is an obvious definition of the composition of morphisms with connection.

Isomorphism classes of gerbes with connection $[\cH,\omega]$ are classified by
the differential cohomology classes $\widehat{DD}[\cH,\omega]\in \widehat{H\Z^3}(M)$, see \cite{MR1876068}, \cite{MR1197353}. In terms of the structure maps $R,I,a$ of differential cohomology we have
$$R(\widehat{DD}[\cH,\omega])=R^\omega\ , \quad I(\widehat{DD}[\cH,\omega])=DD(\cH)\ ,\quad 
\widehat{DD}[\cH,\omega]+a(\alpha)=\widehat{DD}[\cH,\omega+\alpha]\ ,$$

 Following \cite[Def. 1.3]{waldorf1} we make the following definition:
\begin{ddd}
A geometric multiplicative gerbe on a Lie group $G$
is a multiplicative gerbe $(\cG,\mu,a)$  together with a connection
$\omega_\cG$ on $\cG$ and a connection
$\omega_\mu$ on $\mu$ such that
\begin{enumerate}
\item The curvature $\rho$ of $\omega_\mu$ satisfies
\begin{equation}\label{uidqwdwqdqwdqwdqwdqwdqwdqw}
\pr_{23}^*\rho+m_{23}^*\rho=\pr_{12}^*\rho+m_{12}^*\rho\ ,
\end{equation}\ ,
\item $(\mu,\omega_\mu)$ induces an isomorphism of gerbes with connection
\begin{equation}\label{uiwdqwdqwdqwdwd}
(\cG,\omega_\cG)_1\otimes (\cG,\omega_\cG)_2\to (m^*\cG,m^*\omega_\cG+\rho)\ ,
\end{equation}and
\item the associativity $2$-morphism $a$ preserves the connections. 
\end{enumerate}
\end{ddd}

Let $G$ be a compact, connected, and simply connected Lie group. 
For  $x\in H^3(G;\Z)$ there exists a unique bi-invariant form
$\omega_x\in \Omega^3(G)$ which represents the image of $x$ in de Rham cohomology. 
Since $H^2(G;\R/\Z)=0$ and $H^3(G;\Z)$ is torsion-free a class   $\hat x\in \widehat{H\Z}^3(G)$ is completely determined by the invariant $R(\hat x)\in \Omega^3(G)$. Therefore a gerbe $\cG$ on $G$ with $x=DD(\cG)$ has a unique connection $\omega_\cG$
with curvature $R^{\omega_\cG}=\omega_{x}$. 
The following Lemma has been
shown in 
\cite[Sec. 1]{waldorf1}. It extends Lemma \ref{udiqwdqwdwqddqwdqwdqwd} to the geometric context.
\begin{lem}\label{uiwefwefwefewefwefwfewf}
Let $\cG$ be a gerbe on a compact, connected, and simply connected Lie group
with connection $\omega_\cG$ with bi-invariant curvature. This structure can be extended in a unique (up to isomorphism) way to a structure of  a geometric multiplicative gerbe.
\end{lem}
\proof
We will need some of the details of the proof later in the proofs of Lemmas \ref{udqiwdqwdqwdqwdqwdqwd} and \ref{8qqdqwdwqd}.
We must add a connection on the multiplication map (\ref{udiqwdqwdwqddqwdqwdqwd}) which turns $\cG$ into a geometric multiplicative gerbe. For simplicity, we assume  that $G$ is simple.  Then we  have an explicit formula for $\omega_x$ in terms of the Maurer-Cartan form $\theta:=g^{-1}dg\in \Omega^1(G,Lie(G))$. We have $$\omega_x=\frac{k_x}{6c_G}\lk \theta,[\theta,\theta]\rk\ ,$$
where $k_x\in \Z$ depends on $x=DD(\cG)$, $\lk.,.\rk$ is the Killing form,
and $c_G\in \R$ is defined such that
$ \frac{1}{6c_G}\lk \theta,[\theta,\theta]\rk$ represents the image of the generator
of $H^3(G;\Z)\cong \Z$ in de Rham cohomology.
In this case one  can choose  \cite[(1.7)]{waldorf}
\begin{equation}\label{qfwefwefewfewf}
\rho:=\frac{k_x}{c_G}\lk \pr_1^*\theta,\pr_2^*\bar \theta\rk\ ,
\end{equation}
where $\bar \theta:=dgg^{-1}$. By a calculation
one checks that  (\ref{udiqwdqwdwqdwqd}) holds true.

In order to construct the connection on the morphism $\mu$ we use the fact that a again class $\hat y\in \hHZ^3(G^2)$ is completely determined by its
curvature $R(\hat y)\in \Omega^3(G^2)$.
By a calculation we check that
$$\pr_1^*\omega_x+\pr_2^*\omega_x-m^*\omega_x=d\rho\ .$$
Therefore the two sides of the arrow (\ref{uiwdqwdqwdqwdwd})
have the same curvature and hence are isomorphic. Since the $1$-morphism
$\mu$ is the unique one up to isomorphism we can find a unique connection $\omega_\mu$
such that (\ref{uiwdqwdqwdqwdwd}) is an isomorphism of gerbes with connection.

The relation (\ref{uidqwdwqdqwdqwdqwdqwdqwdqw}) implies that
 the bundle 
$\frac{t_1}{t_0}$ is flat, where $t_0,t_1$ are as in (\ref{zuiqwkdqwdqwdqwdqwdd}). Since $G^3$ is simply-connected 
the bundle is trivial. If in the proof of Lemma \ref{udiqwdqwdwqddqwdqwdqwd} we take  a trivialisation $\tilde a$ which preserves the connection, then
its deviation from satisfying the higher coherence relation is a constant group  cocycle $c_{\tilde a}$ with values in $U(1)$, in other words a cocycle in
$C_{gr}^4(G,U(1)^\delta)$, where $U(1)^\delta$ has the discrete topology.
Since $G$ is connected  we have (as for the coefficients $\Z$) for $i\ge 1$ that
$H^i_{gr}(G;U(1)^\delta)=H^i_{gr}(\{1\};U(1)^\delta)=0$.
Therefore we can get an associativity morphism $a:=\tilde a d^{-1}$ for a unique constant $d\in C_{gr}^3(G,U(1)^\delta)\cong U(1)$ such that $\delta d=c_{\tilde a}$. Note that $a$ still preserves the connection.
\hB

Let us apply this to the basic gerbe $\cG$ on $Spin(n)$ which by Lemma \ref{uiwefwefwefewefwefwfewf} becomes
a geometric multiplicative gerbe. Its bi-invariant  curvature will be denoted by 
$CS\in \Omega^3(Spin(n))$.

\subsection{Geometric string structures}\label{str111}

 First we recall the notion of a string structure  according to \cite{waldorf}.
We consider an $n$-dimensional spin vector bundle $V\to M$.  Let $p:P\to M$ be the corresponding $Spin(n)$-principal bundle   (earlier we used the longer notation $P=Spin(V)$). Then we have a canonical isomorphism
\begin{equation}\label{tzzttzuiiu}
(\id,g):P\times_M P\stackrel{\sim}{\to} P\times Spin(n)
\end{equation} of $Spin(n)$-principal bundles whose inverse is given by 
$(p,h)\mapsto (p,ph)$.

Let $\cG$ be the basic multiplicative gerbe on $Spin(n)$.
We define the gerbe
$\cP:=g^*\cG$ on $P\times_MP$. The multiplicative structure of $\cG$ induces
a $1$-morphism 
\begin{equation}\label{uqidqwdqwdwqdwdqdwq}
\nu: \cP_{12}\otimes \cP_{23}\to \cP_{13}
\end{equation}
together with an associativity $2$-morphism, where $\pr_{ij}:P\times_MP\times_MP\to P\times_MP$ are the projections and $\cP_{ij}:=\pr_{ij}^*\cP$.
In detail, if we define $g_{ij}:=g\circ \pr_{ij}$, then we have
$m\circ (g_{12}\times g_{23})=g_{13}$ and $\cP_{ij}\cong g_{ij}^*\cG$ so that
$P_{13}\cong (g_{12}\times g_{23})^*m^*\cG$ and $P_{12}\otimes \cP_{23}\cong (g_{12}\times g_{23})^*(\cG_1\otimes \cG_2)$.
Therefore $\nu$ is defined as the pull-back of $\mu$ in (\ref{iudqwdwqdwqd})
via $g_{12}\times g_{23}$.

According to \cite[Sec. 2.1]{waldorf} we make the following definition:
\begin{ddd}
The gerbe $\cP$ over $P$ together with the multiplication $1$-morphism (\ref{uqidqwdqwdwqdwdqdwq}) and the associativity $2$-morphism
 is called   the Chern-Simons bundle $2$-gerbe  $\cCS$ on $M$ associated to the spin bundle $V$.
\end{ddd}

By \cite[Thm 1.1.4]{waldorf} we can
define a string structure as a trivialisation of the Chern-Simons bundle $2$-gerbe. In detail, for $i=1,2$ let $\pr_i:P\times_MP\to P$ denote the projections and set $\cS_i:=\pr_i^*\cS$.
\begin{ddd} A string structure on the spin bundle $V$
is a gerbe $\cS$ over $P$ together with a $1$-morphism
$$f:\cP\otimes \cS_2\to \cS_1$$ and an  associativity $2$-morphism (which essentially turns $\cS$  into module over $\cP$) satisfying a higher coherence condition \cite[Def. 2.2.1]{waldorf}.
\end{ddd}
We will usually denote a string structure by the same symbol as its underlying gerbe.
By \cite[Lemma 2.2.2]{waldorf} a string structure exists if and only if $\frac{p_1}{2}(V)=0$. In this case the isomorphism classes of string structures form a torsor over $H^3(M;\Z)$, see  \cite[Thm.1.1.2]{waldorf} where this result is attributed to \cite{MR2032513}.
We continue with recalling notions and results from \cite{waldorf}. 
\begin{ddd}\label{uiufwefwfef4}
A connection $h$  on the bundle $2$-gerbe $\cCS$ 
consists of
\begin{enumerate}
\item a  $3$-form $\kappa_h\in \Omega^3(P)$, 
\item a connection $\omega_h$ on the gerbe $\cP$,
\item a connection $\sigma_h$ on the multiplication (\ref{uqidqwdqwdwqdwdqdwq})
\end{enumerate}  
such that
\begin{enumerate}
\item\label{zufuewfeiooiwoef}  $\pr_2^*\kappa_h-\pr_1^*\kappa_h=R^{\omega_h}$,
\item\label{z8idqwdqwdwqdw} $(\nu,\sigma_h)$ realizes an isomorphism of gerbes with connection, and
\item\label{uiudqwdqwdqwdwqdwqdwd}  the associativity morphism preserves connections.
\end{enumerate} 
\end{ddd}
It follows that $\pr_1^* d\kappa_h=\pr_2^* d\kappa_h$. Therefore there exists  a unique closed $4$-form
$R^h\in \Omega^4(M)$ such that $p^*R^h=d\kappa_h$. The cohomology class of 
$R^h$ is the image of $\frac{p_1}{2}(V)$ in de Rham cohomology.

\begin{lem}[\cite{waldorf1}]\label{udqiwdqwdqwdqwdqwdqwd} Let   $\bV=(V,\nabla^V,h^V)$ be a geometric spin bundle.
Then we have an associated  connection $h_\bV$ on the Chern-Simons $2$-gerbe $\cCS$ of $V$.
\end{lem}
\proof
In the following we describe the construction of $h_\bV$ which is due to \cite[Sec. 3]{waldorf1}. We will need the details later in the proof of Lemma \ref{8qqdqwdwqd}.
Recall that the basic gerbe $\cG$ on $Spin(n)$ has a unique connection
$\omega_\cG$ with curvature 
$$R^{\omega_\cG}=CS=\frac{1}{6c_{Spin(n)}}\lk\theta,[\theta,\theta]\rk\in \Omega^3(Spin(n))\ .$$ 
The gerbe  $\cP=g^*\cG$ has an induced connection
$\omega_\cP=g^*\omega_\cG$ with curvature $g^*CS$.

The bundle $p^*P\cong P\times_MP$ has a canonical trivialisation (\ref{tzzttzuiiu})
which we  denote by $Q$ for the moment. It induces a taming $(p^*\cW)_{t_Q}$ as explained in Section \ref{ll}. We have seen in Lemma \ref{uqifqffff77wefefwf} that  $CS(\nabla^{p^*V})=\frac{1}{2}\eta^3((p^*\cW)_{t_Q})\in \Omega^3(P)$ is the usual Chern-Simons form of $\nabla^{p^*V}$ (in the trivialisation given by $Q$).
Let $A\in \Omega^1(P,spin(n))$ denote the connection one-form of $\nabla^V$.
Then in the new notation
$$CS(\nabla^{p^*V})=\frac{1}{c_{Spin(n)}} \left(\lk dA,A \rk+\frac{2}{3} \lk A,\frac{1}{2}[A,A]\rk\right)\ ,$$
and this fixes the sign of $c_{Spin(n)}$. 
Then we define 
$$\omega:=-\frac{1}{c_{Spin(n)}}\lk \pr_1^*A,g^*\bar \theta \rk\in  \Omega^2(P\times_MP)\ .$$
We have  by a direct calculation (see  \cite[Sec. 3.1]{waldorf})
\begin{equation}\label{zwdwqdwqdwqdwqd}
-\pr_2^*CS(\nabla^{p^*V})+\pr_1^*CS(\nabla^{p^*V})=g^*CS+d\omega\ .
\end{equation}
We set $\kappa_{h_{\bV}}:=-CS(\nabla^{p^*V})$ and $\omega_{h_\bV}:=\omega_{\cP}+\omega$. Then condition \ref{zufuewfeiooiwoef} holds true.
Since the multiplication $\nu$ in (\ref{uqidqwdqwdwqdwdqdwq}) is defined as the pull-back of the multiplicative structure $\mu$ of $\cG$ we can define the connection $\sigma_{h_{\bV}}$ on $\nu$ by pulling back the connection $\omega_\mu$.
Since the associativity morphism for the Chern-Simons gerbe is also obtained by pulling back the associativity morphism of the multiplicative structure on $\cG$
our definition of $\sigma_{h_{\bV}}$ ensures that the associativity
morphism of the Chern-Simons gerbe preserves connections, hence condition \ref{uiudqwdqwdqwdwqdwqdwd} holds true.
Moreover it satisfies the part \ref{uiduwqidhqwdwqdwqdqd} of  condition Definition \ref{uiufwefwfef4}.\ref{z8idqwdqwdwqdw}.
In order to verify condition \ref{uiufwefwfef4}.\ref{z8idqwdqwdwqdw} completely we must check
that $(\nu,\sigma_{h_{\bV}})$ is flat. With the curvature $\rho$ of $\mu$ given by (\ref{qfwefwefewfewf}) this is the equality 
$$(g_{12}\times g_{23})^*\rho=\pr_{12}^*\omega+\pr_{23}^*\omega-\pr_{13}^*\omega$$
which again can be checked by a direct calculation (compare \cite[(3.21)]{waldorf1})
\hB

We now consider a Chern-Simons gerbe $\cCS$ with a connection $h$. 
\begin{ddd}
A geometric string structure is a triple $str:=(\cS,\omega_\cS,\omega_f)$
of a string structure $\cS$ with action $f:\cP\otimes \cS_2\to \cS_1$
together with a connection $\omega_\cS$ on $\cS$ and a connection $\omega_f$ on the morphism $f$ such that
\begin{enumerate}
\item $(f,\omega_f)$ realizes an isomorphism of gerbes with connection, and
\item the associativity $2$-morphism  preserves connections.
\end{enumerate}
\end{ddd}

It is shown in  \cite[Thm 1.3.4]{waldorf} that a string structure $\cS$ can always be refined to a
geometric string structure.

Assume that we have chosen a  geometric string structure $str$. It was shown in \cite[Thm 1.3.3]{waldorf} that there is a unique form $H_{str}\in \Omega^3(M)$ such that
\begin{equation}\label{zdqwdqwddddqwdq}
p^*H_{str}=R^{\omega_{\cS}}+\kappa_h\ .
\end{equation}
This form is closely related to the Cheeger-Simons cocycle
$ \frac{\hat p_1}{2}(\bV):Z_3(E)\to U(1)$. Indeed, if $z\in Z_3(E)$, then there exists a neighbourhood $U\subseteq E$ of  the trace $|z|$ of $z$  such that 
that $\frac{p_1}{2}(V)_{|U}=0$. Therefore there exists a geometric string structure
$str$ on $P_{U}$. Then 
$$ \frac{\hat p_1}{2}(\bV)(z)=\exp(2\pi i \int_z H_{str})\ .$$
 This follows by combining
(\ref{uwidqwuidhqwdqwdqdqd}) with 
equation (\ref{u8wfe9fwefewfwf}) below (compare \cite[Sec. 3.4]{waldorf}).

\subsection{The string structure associated to a trivialisation}

In this subsection we show that a trivialisation $Q$ of the geometric spin bundle $\bV$ on $M$ induces a geometric string structure $str_Q$.  In Lemma \ref{8qqdqwdwqd}  we calculate the associated form $H_{str_Q}\in \Omega^3(M)$ (see (\ref{zdqwdqwddddqwdq})).

\begin{lem}\label{zugdqwuzdwdwqd}
A trivialisation $Q:P\to M\times Spin(n)$ of $Spin(n)$-principal bundles gives rise to a string structure $\cS_Q$.
\end{lem}
\proof
We first consider the trivial bundle $\R^n\to *$ with spin structure $Spin(n)\to *$.
Its Chern-Simons gerbe is $\tilde m^*\cG$, where $\tilde m:Spin(n)\times Spin(n)\to Spin(n)$ is given by $\tilde m(g,h):=g^{-1}h$.
The multiplicative gerbe $\cG$ on $Spin(n)$ can be considered as a 
string structure.  

The trivialization $Q$ of $P$ gives a pull-back diagramm
$$\xymatrix{P\ar[d]^p\ar[r]^q&Spin(n)\ar[d]\\M\ar[r]&{*}}$$
of $Spin(n)$-principal bundles. The Chern-Simons bundle $2$-gerbe of $V$ as well as the string structure are now obtained by induced pull-backs.

More explicitly, we have $\cP\cong (q_1\times q_2)^*\tilde m^*\cG$ and set $\cS_Q:=q^{-1*}\cG$. The action map is then given by
\begin{equation}\label{zuqdqwdqwdqwdd}
f:=(\tilde m\times q_2^{-1})^*\mu\ ,
\end{equation}
where we use the identity $q_1^{-1}=m\circ (\tilde m\circ (q_1\times q_2)\times  q_2^{-1})$.
\hB

 Let $\bV$ be a geometric spin bundle and $h_\bV$ be the associated connection on $\cCS$.
 In Lemma \ref{zugdqwuzdwdwqd}
 we have seen that a trivialisation $Q$ of $P$ gives a string structure $\cS_Q=q^*\cG$, where $q:P\stackrel{\sim}{\to} M\times Spin(n)\to Spin(n)$. 

\begin{lem} \label{8qqdqwdwqd}
A trivialisation $Q$ gives a natural extension of $\cS_Q$ to a geometric string structure $str_Q$. 
Moreover, we have 
\begin{equation}\label{u8wfe9fwefewfwf}
H_{str_Q}= -\frac{1}{2} \eta^3(\cW_{t_Q})\ .
\end{equation}
\end{lem}
\proof
Using the trivialization $Q$ we identify (with $G:=Spin(n)$)
$$P\cong M\times G\ ,\quad P\times_MP\cong M\times G\times G$$
so that
$$q(b,h)=h\ ,\quad \pr_1(b,h,l)=(b,h)\ ,\quad \pr_2(b,h,l)=(b,l)\ , \quad g(b,h,l)= \tilde m(h,l)=h^{-1}l\ .$$
We define the form
$$\alpha:=\frac{1}{c_{Spin(n)}}\lk A,q^*\theta \rk\in \Omega^2(P)\ .$$
We equip the gerbe $\cS_Q$ with the connection
$$\omega_{S_Q}:=q^{-1*}\omega_\cG+\alpha\ .$$ Furthermore, the action 
(\ref{zuqdqwdqwdqwdd}) will be equipped with the connection
$\omega_f:=(\tilde m\times q_2^{-1})^*\omega_\mu$.
This already ensures that the associativity morphism obtained by pull-back from
 the multiplicative structure of $\cG$ respects the connection. Furthermore part \ref{uiduwqidhqwdwqdwqdqd} of the compatibility of $\omega_f$ is satisfied.
It remains to show that
$(f,\omega_f)$ is flat, i.e. we must check the identity of $2$-forms on $P\times_MP$:
$$(\tilde m\times q_2^{-1})^*\rho-\omega-\alpha_2+\alpha_1=0\ ,$$
where we use the notation $\alpha_i=\pr_i^*\alpha$.
This is a straightforward  calculation.

We now show (\ref{u8wfe9fwefewfwf}).
Note that we have $R^{\omega_{S_Q}}=q^*CS+d\alpha$.
 Let $s_Q:M\to P$ denote the section given by $Q$.
Then the composition $q\circ s_Q:M\to Spin(n)$ is constant and hence
$s_Q^* q^*CS=0$ and $s_Q^*\alpha=\frac{1}{c_{Spin(n)}}\lk s_Q^*A,s_Q^*q^*\theta \rk =0$.
Moreover, the pull-back via $s_Q$ of the canonical trivialisation of $p^*V$ is exactly the trivialisation of $V$ given by $Q$.
Therefore
$s_Q^*(CS(\nabla^{p^*V}))=CS(\nabla^V)=\frac{1}{2}\eta^3(\cW_{t_Q})$, the Chern-Simons form of the connection $\nabla^V$ in the trivialisation $Q$.
Applying $s^*_Q$ to $p^*H_{str_Q}=q^*CS-CS(\nabla^{p^*V})$ we get
 $H_{str_Q}=-\frac{1}{2}\eta^3(\cW_{t_Q})$.
\hB

\subsection{Proof of Theorem  \ref{tzu}}

Let $\pi:E\to B$ be our surface bundle with 
the geometric   spin bundle $\bV$  on $E$.
Let $\bL$  be the geometric line bundle on $B$ constructed 
in Subsection  \ref{ztgudqwdqwd}.
\begin{theorem}\label{zuadddqwdqd}
 A geometric string structure $str=(\cS,\omega_\cS,\omega_f)$ on $\bV$ gives a functorial unit-norm section $s_{str}\in C^\infty(B,L)$.  It satisfies
$$\nabla^L \log s_{str}=2\pi i \int_{E/B}H_{str}\ .$$
\end{theorem}
The meaning of the adjective {\em functorial} is here again  the obvious compatibility of the construction
with cartesian diagrams of the form (\ref{uqidqwdqwdwqdqd545454545}).

\proof
If $U\subseteq B$ is an contractible open subset and $Q\in I_U$ is a trivialisation of $P_{E_U}$ , then by the construction of $L$ we have a section
$s_Q\in C^\infty(U,L)$.
Using the string structure we will define a function $a_Q\in C^\infty(U,U(1))$ such that
$\tilde s_Q:=a_Q s_Q$ is independent of the choice of $Q$.
The collection $(\tilde s_Q)_{U\subseteq B,Q\in I_U}$ therefore
defines a global section $s_{str}\in C^\infty(B,L)$. In order to show the second part we calculate that
$$\nabla^L \log \tilde s_Q=2\pi i\int_{E_U/U} H_{str}\ .$$

In order to define $a_Q$ we consider the projection
$\Pr:[0,1]\times E_U
\to E_U$. There exists a unique string structure $\tilde \cS$ on $\Pr^* V$
which restricts to $\cS_{E_U}$ on $\{1\}\times E_U$, and to $\cS_Q$ on $\{0\}\times E_U$.
In fact, since $H^3(E_U;\Z)=0=H^3([0,1]\times E_U;\Z)$ there is  only one upto  isomorphism string structure on $V_{E_U}$ and $\Pr^*V_{E_U}$, respectively.
We first fix the isomorphisms above. 
Then we can choose a geometric string structure
$\widetilde{str}=(\tilde \cS,\tilde \omega_{\tilde \cS},\tilde \omega_{\tilde f})$
which restricts to the given ones $str_Q$ (defined in Lemma \ref{8qqdqwdwqd}) and $str_{|E_U}$ at the endpoints $\{0\}\times E_U$ and $\{1\}\times E_U$. This follows from the fact
that different geometric refinements of a string structure can be glued using a partition of unity (use \cite[Prop. 3.3.4]{waldorf}).

We define
$$a_Q:=\exp\left(2\pi i \int_{[0,1]\times E_U/U} H_{\widetilde{str}}\right)\ ,$$ 
where $H_{\widetilde{str}}$ is the $3$-form associated to the string structure$ \widetilde{str}$ by 
(\ref{zdqwdqwddddqwdq}).

We first show that $a_Q$ does not depend on the choices made in the construction, namely the isomorphisms of gerbes $\tilde \cS_{\{0\}\times E_U}\cong \cS_Q$, $\tilde \cS_{\{1\}\times E_U}\cong \cS_{E_U}$,
and of the geometry. From two such choices $str,str^\prime$ we can produce a geometric string structure $\widehat{str}$
on $$\pr_{E_U}^*V\to S^1\times E_U\cong ([0,1]\times E_U)\sqcup ([0,1]\times E_U)/\sim
\ ,$$where $\sim$ identifies the endpoints, such that
$$\int_{[0,1]\times E_U/U} H_{\widetilde{str}}-\int_{[0,1]\times E_U/U} H_{\widetilde{str}^\prime}=\int_{S^1\times E_U/U} H_{\widehat{str}}\ .$$
Over each point $u\in U$ the right-hand side $$\exp\left(2\pi i \int_{S^1\times E_U/U} H_{\widehat{str}}\right)$$ is the evaluation
$ \frac{\hat p_1}{2} (\pr_{E_U}^*\bV)(z_u)\in U(1)$ of the Cheeger-Simons character of $ \frac{\hat p_1}{2} (\pr_{E_U}^*\bV)$ on the cycle $$z_u=(S^1\times E_{\{u\}}\to D^2\times  E)\in Z_3(D^2\times E)\ .$$
Note that $z_u$ is the boundary of the $4$-chain $(\phi:D^2\times E_{\{u\}}\to E)\in C_4(D^2\times E)$, and therefore
$$\frac{\hat p_1}{2}(\pr_{E_U}^*\bV)(z_u)=\exp\left(\pi i\int_{D^2\times E_{\{u\}}} \phi^* p_1(\nabla^V)\right)=1\ .$$
This implies $$\int_{[0,1]\times E_U/U} H_{\widetilde{str}}-\int_{[0,1]\times E_U/U} H_{\widetilde{str}^\prime}\in C^\infty(U,\Z)\ .$$
 
Next we calculate the quotient
$\frac{a_{Q^\prime}}{a_Q}\in C^\infty(U,U(1))$.
We choose a homotopy $H$ from $Q$ to $Q^\prime$.
This homotopy gives a geometric string structure
$str_H$ on $[0,1]\times E_U$ which connects
$str_Q$ and $str_{Q^\prime}$.
We consider the projection 
$$\Pr_{E_U}:[0,1]\times [0,1]\times E_U\to E_U\ $$
On $ \Pr_{E_U}^*V$  we can find a geometric string structure
$\widehat{str}$ which restricts to 
$\widetilde{str}$ on $\{0\}\times [0,1]\times E_U$, to
$\widetilde {str}^\prime$ on $\{1\}\times [0,1]\times E_U$, to $str_H$  in $[0,1]\times \{0\}\times E_U$,  and to
$\pr_{E_U}^* str_{E_U}$ on $[0,1]\times \{1\}\times E_U$.
By Stoke's theorem
\begin{eqnarray*}\lefteqn{-
\int_{[0,1]\times E_U/U} H_{\widetilde{str}}+\int_{[0,1]\times E_U/U} H_{\widetilde{str}^\prime}+\int_{[0,1]\times E_U/U} H_{str_H}}&&\\&=&\int_{[0,1]\times [0,1]\times E_U/U} dH_{\widetilde{str}}\\
&=& \frac{1}{2}\int_{[0,1]\times [0,1]\times E_U/U}  \Pr_{E_U}^* p_1(\nabla^V)\\
&=&0\ .
\end{eqnarray*} 
This implies that
$$\frac{a_{Q^\prime}}{a_{Q}}=\exp\left( -2 \pi i \int_{[0,1]\times E_U/U} H_{str_H}\right)=\exp\left( \pi i \int_{[0,1]\times E_U/U}  \eta^3(\Pr^* \cW_{t_H})\right)=\frac{1}{c(Q^\prime,Q)}\ ,$$
where $c(Q^\prime,Q)\in C^\infty(U,U(1))=\frac{s_Q^\prime}{s_Q}$ is as in (\ref{uiwqddqwd44}). We thus get
$$\frac{\tilde s_Q^\prime}{\tilde s_Q}=c(Q^\prime,Q) \frac{a_{Q^\prime}}{a_{Q}}=1$$
as required.
We now calculate the covariant derivative of $\tilde s_Q$.
We use (\ref{uiwqddqwd444}) and Stoke's theorem
\begin{eqnarray*}
\nabla^L \log \tilde s_Q&=&\nabla^L \log s_Q+2\pi i d\int_{[0,1]\times E_U/U} H_{\widetilde{str}}\\&=&-\pi i \int_{[0,1]\times E_U/U}  \eta^3(\cW_{t_Q})+2\pi i \int_{E_U/U} H_{str}- 2\pi i\int_{E_U/U} H_{str_Q}\\&=&
2\pi i \int_{E_U/U} H_{str}
\end{eqnarray*}
\hB


\begin{thebibliography}{CJM{\etalchar{+}}05}

\bibitem[AS68]{MR0236950}
M.~F. Atiyah and I.~M. Singer.
\newblock The index of elliptic operators. {I}.
\newblock {\em Ann. of Math. (2)}, 87:484--530, 1968.

\bibitem[BC89]{MR966608}
J.-M. Bismut and J. Cheeger.
\newblock {$\eta$}-invariants and their adiabatic limits.
\newblock {\em J. Amer. Math. Soc.}, 2(1):33--70, 1989.

\bibitem[BGV92]{MR1215720}
N.~Berline, E.~Getzler, and M.~Vergne.
\newblock {\em Heat kernels and {D}irac operators}, volume 298 of {\em
  Grundlehren der Mathematischen Wissenschaften [Fundamental Principles of
  Mathematical Sciences]}.
\newblock Springer-Verlag, Berlin, 1992.

\bibitem[BKS09]{bks}
U.~ Bunke, M.~Kreck, and Th.~Schick.
\newblock A geometric description of smooth cohomology, 2009, To appear.
 arXiv.org:0903.5290.

\bibitem[Bor92]{MR1186039}
D.~Borthwick.
\newblock The {P}faffian line bundle.
\newblock {\em Comm. Math. Phys.}, 149(3):463--493, 1992.

\bibitem[Bry93]{MR1197353}
J.~L. Brylinski.
\newblock {\em Loop spaces, characteristic classes and geometric quantization},
  volume 107 of {\em Progress in Mathematics}.
\newblock Birkh\"auser Boston Inc., Boston, MA, 1993.

\bibitem[BS07]{bunke-20071}
U.~Bunke and Th.~Schick.
\newblock Smooth K-theory, 2007, To appear.
arXiv.org:0707.0046.

\bibitem[BS08]{MR2443109}
U.~Bunke and Th.~Schick.
\newblock Real secondary index theory.
\newblock {\em Algebr. Geom. Topol.}, 8(2):1093--1139, 2008.

\bibitem[BS09]{bunke-2009}
U.~Bunke and Th.~Schick.
\newblock Uniqueness of smooth extensions of generalized cohomology theories,
  2009, Submitted. arXiv.org:0901.4423.

\bibitem[BSST08]{bunke-2008-2006}U.~Bunke, Th.~Schick, M.~Spitzweck, and A.~Thom.
\newblock {\em Duality for topological abelian group stacks and T-duality}.
\newblock Cortinas, Guillermo (ed.) et al., $K$-theory and noncommutative
geometry. Proceedings of the ICM 2006 satellite conference, Valladolid, Spain,  2006. 
\newblock Z\"urich: European Mathematical Society
(EMS). Series of Congress Reports, 227-347 (2008). 
 arXiv.org:math/0701428.

\bibitem[Bun]{bunke-chern}
U.~Bunke.
\newblock Chern classes on differential K-theory, 2009, Submitted.
   arXiv.org:0907.2504.

\bibitem[Bun09]{MR2191484}
U.~Bunke.
\newblock Index theory, eta forms, and {D}eligne cohomology.
\newblock {\em Mem. Amer. Math. Soc.}, 198(928):vi+120, 2009.

\bibitem[CJM{\etalchar{+}}05]{MR2174418}
A.~L. Carey, St.~Johnson, M.~K. Murray, D.~Stevenson, and B.~L.
  Wang.
\newblock Bundle gerbes for {C}hern-{S}imons and {W}ess-{Z}umino-{W}itten
  theories.
\newblock {\em Comm. Math. Phys.}, 259(3):577--613, 2005.

\bibitem[CS85]{MR827262}
J.~Cheeger and J.~Simons.
\newblock Differential characters and geometric invariants.
\newblock In {\em Geometry and topology (College Park, Md., 1983/84)}, volume
  1167 of {\em Lecture Notes in Math.}, pages 50--80. Springer, Berlin, 1985.

\bibitem[Del87]{MR902592}
P.~Deligne.
\newblock Le d\'eterminant de la cohomologie.
\newblock In {\em Current trends in arithmetical algebraic geometry ({A}rcata,
  {C}alif., 1985)}, volume~67 of {\em Contemp. Math.}, pages 93--177. Amer.
  Math. Soc., Providence, RI, 1987.
\bibitem[Fra91]{MR1085257}
J.~Franke.
\newblock Chern functors.
\newblock In {\em Arithmetic algebraic geometry ({T}exel, 1989)}, volume~89 of
  {\em Progr. Math.}, pages 75--152. Birkh\"auser Boston, Boston, MA, 1991.




\bibitem[FL09]{freed-2009}
D.~S. Freed and J.~Lott.
\newblock An index theorem in differential K-theory, 2009,
 arXiv.org:0907.3508.

\bibitem[FM06]{MR2207325}
D.~S. Freed and G.~W. Moore.
\newblock Setting the quantum integrand of {M}-theory.
\newblock {\em Comm. Math. Phys.}, 263(1):89--132, 2006.

\bibitem[Fre03]{MR915823}
D.~S. Freed.
\newblock On determinant line bundles.
\newblock In: {\em Mathematical aspects of string theory ({S}an {D}iego,
              {C}alif., 1986)}, pages 189--238, World Sci. Publishing, Singapore 1987.
 

\bibitem[Fre02]{MR1971377}
D.~S. Freed.
\newblock {$K$}-theory in quantum field theory.
\newblock In {\em Current developments in mathematics, 2001}, pages 41--87.
  Int. Press, Somerville, MA, 2002.

\bibitem[Hei05]{MR2206877}
J.~Heinloth.
\newblock Notes on differentiable stacks.
\newblock In {\em Mathematisches {I}nstitut, {G}eorg-{A}ugust-{U}niversit\"at
  {G}\"ottingen: {S}eminars {W}inter {T}erm 2004/2005}, pages 1--32.
  Universit\"atsdrucke G\"ottingen, G\"ottingen, 2005.

\bibitem[Hit01]{MR1876068}
N.~Hitchin.
\newblock Lectures on special {L}agrangian submanifolds.
\newblock In {\em Winter {S}chool on {M}irror {S}ymmetry, {V}ector {B}undles
  and {L}agrangian {S}ubmanifolds ({C}ambridge, {MA}, 1999)}, volume~23 of {\em
  AMS/IP Stud. Adv. Math.}, pages 151--182. Amer. Math. Soc., Providence, RI,
  2001.

\bibitem[HS05]{MR2192936}
M.~J. Hopkins and I.~M. Singer.
\newblock Quadratic functions in geometry, topology, and {M}-theory.
\newblock {\em J. Differential Geom.}, 70(3):329--452, 2005.

\bibitem[Mad09]{madsen}
I.~Madsen.
\newblock An integral Riemann-Roch theorem for surface bundles, 2009,
 arXix.org:0901.4240.

\bibitem[SS08]{MR2365651}
J.~Simons and D.~Sullivan.
\newblock Axiomatic characterization of ordinary differential cohomology.
\newblock {\em J. Topol.}, 1(1):45--56, 2008.

\bibitem[Ste04]{MR2032513}
D.~Stevenson.
\newblock Bundle 2-gerbes.
\newblock {\em Proc. London Math. Soc. (3)}, 88(2):405--435, 2004.

\bibitem[Wala]{waldorf1}
K.~Waldorf.
\newblock Multiplicative bundle gerbes with connection,
 arXix.org:0804.4835v3.

\bibitem[Walb]{waldorf}
K.~Waldorf.
\newblock String connections and Chern-Cimons theory,
 arXix.org:0906.0117.

\bibitem[Wal07]{MR2318389}
K.~Waldorf.
\newblock More morphisms between bundle gerbes.
\newblock {\em Theory Appl. Categ.}, 18:No. 9, 240--273 (electronic), 2007.

\end{thebibliography}
 \newcommand{\etalchar}[1]{$^{#1}$}

\end{document}